\documentclass[11pt,reqno]{amsart}
\usepackage{t1enc}
\usepackage[latin1]{inputenc}
\usepackage[english]{babel}
\usepackage{amsmath,amsthm}
\usepackage{amsfonts}
\usepackage{latexsym}
\usepackage[dvips]{graphicx}
\usepackage{graphicx}
\usepackage{float}
\usepackage[natural]{xcolor}
\usepackage{algorithm}
\usepackage{algorithmic}
\usepackage{enumerate}
\usepackage{appendix}
\usepackage{multirow}
\usepackage{xcolor}
\usepackage[colorlinks,linkcolor=blue]{hyperref}
\usepackage{setspace}
\usepackage{comment}
\usepackage{fullpage}
\usepackage[shortlabels]{enumitem}

\usepackage{listings}
\newtheorem{theorem}{Theorem}\numberwithin{theorem}{section}
\newtheorem{lemma}[theorem]{Lemma}\numberwithin{theorem}{section}
\newtheorem{proposition}[theorem]{Proposition}
\numberwithin{theorem}{section}

\newtheorem{definition}[theorem]{Definition}

\newtheorem{claim}[theorem]{Claim}

\def\int{\textrm{int}}

\begin{document}

\onehalfspacing
\title{Perturbation of dense graphs}

\author{Jie Han}
\address{JH, BW and JZ. School of Mathematics and Statistics, Beijing Institute of Technology, China\\
Email: \texttt{(JH) han.jie@bit.edu.cn, (BW) bin.wang@bit.edu.cn, (JZ) jxuezhang@163.com.}}
\author{Seonghyuk Im}
\address{SI. Department of Mathematical Science, KAIST and Extremal Combinatorics and Probability Group (ECOPRO), Institute for Basic Science (IBS), South Korea, Email: \texttt{seonghyuk@kaist.ac.kr}.}
\author{Bin Wang}
\author{Junxue Zhang}

\maketitle
\begin{abstract}

In the past two decades, various properties of randomly perturbed/augmented (hyper)graphs have been intensively studied, since the model was introduced by Bohman, Frieze and Martin in 2003.
The model usually considers a deterministic graph $G$ with minimum degree condition, perturbed/augmented by a binomial random graph $G(n,p)$ on the same vertex set.
In this paper, we show that for many problems of finding spanning subgraphs, one can indeed relax the minimum degree condition to a density condition.
This includes the embedding problem for $F$-factors when $F$ is not a forest, graphs with bounded maximum degree, $r$-th power of $k$-uniform tight Hamilton cycles for $r,k\ge 2$, and $k$-uniform Hamilton $\ell$-cycles for $\ell\in[2,k-1]$.
These results strengthen the results of Balogh, Treglown, and Wagner, of B\"{o}ttcher, Montgomery, Parczyk, and Person, and of Chang, Han and Thoma.

\end{abstract}

\section{Introduction}

The binomial random  $k$-graph, since it was introduced  by Erd\H{o}s and R\'{e}nyi
\cite{MR125031} in the 60s, has found applications in various fields of science, such as statstical physics, computer science, and complex networks.
For $k\geq2$, let $G^{(k)}(n,p)$ be the binomial random $k$-graph on a finite set $V$ of $n$ vertices and each $k$-set forms an edge independently with probability $p$ (when $k=2$, we usually abbreviate it as $G(n,p)$).
We say that $G^{(k)}(n,p)$ has some property~$\mathcal{P}$ \emph{with high probability} (w.h.p.) if $\lim_{n\rightarrow\infty}\mathbb{P}[G^{(k)}(n,p)\in \mathcal{P}] = 1$.
In this paper we are interested in various spanning subgraphs of a dense $k$-graph, perturbed by a binomial random $k$-graph on the same vertex set.

Let $k\geq2$, a $k$-graph $G=(V,E)$ consists of a vertex set of order $n$ and an edge set $E\subseteq\binom{V}{k}$.
If $E=\binom{V}{k}$, then we call $G$ a \emph{complete $k$-graph}, denoted by $K_{n}^{(k)}$.
For any $d$-subset $S$ of $V$ where $d\in[k-1]$, the \emph{degree} of $S$, denoted by $\deg_H(S)$ is the number of edges containing $S$.
The \emph{minimum} $d$-\emph{degree} $\delta_d(H)$ is the minimum of $\deg_H(S)$ over all $d$-subsets $S$ of $V$.

\subsection{Thresholds in $G(n,p)$}
Thresholds of graph properties have been one of the central objects in random graph theory.
We say that a function $\hat{p}: \mathbb{N} \rightarrow [0,1]$ is a \emph{threshold} for a graph property~$\mathcal{P}$ if
\begin{equation}
\nonumber
\mathbb{P}[G^{(k)}(n,p)\in \mathcal{P}]=
\left\{
\begin{aligned}
0,\ {\rm if\ } p&/\hat{p}\to 0, \\
1,\ {\rm if\ }p&/\hat{p}\to \infty.
\end{aligned}
 \right.
\end{equation}

A graph is Hamiltonian if it has a cycle that passes through every vertex exactly once.
The threshold for Hamiltonicity was famously determined by P\'{o}sa~\cite{MR389666} and independently by Kor\v{s}unov~\cite{MR434878} in the 70s as $\log n/n$, that is, they showed that $G(n,p)$ contains a Hamilton cycle with high probability if $p \ge \omega(\log n / n)$.
Since then, determining the thresholds for various properties has been a central theme, inspiring a number of problems and results in probabilistic combinatorics.

\subsection{Randomly perturbed/augmented $k$-graphs}\label{subsec:randomly_perturbed_graphs}

A classical result of Dirac~\cite{Dirac} says that every $n$-vertex graph with minimum degree at least $n/2$ is Hamiltonian, given that $n\ge 3$.
The assumption on the minimum degree is best possible, as seen by slightly unbalanced complete bipartite graphs.
However, this extremal graph is highly structured and very far from a typical graph.

Bohman, Frieze, and Martin~\cite{MR1943857} were the first to consider the model of combining a (dense) deterministic graph with a set of random edges.
More precisely, let $G$ be a $k$-graph on $n$ vertices ($k\geq2$), and we consider the union of $G$ and a random $k$-graph $G^{(k)}(n,p)$ defined on the same vertex set $V(G)$.
This is essentially equivalent to adding a uniformly random set of edges to $G$.
Then the resulting graph $G \cup G^{(k)}(n,p)$ is called a \emph{randomly perturbed $k$-graph}.
In particular, Bohman, Frieze, and Martin~\cite{MR1943857} showed that for every $\alpha>0$, there exists $C>0$, such that if an $n$-vertex graph $G$ satisfies $\delta(G)\ge \alpha n$, then w.h.p.~$G \cup G(n, C/n)$ is Hamiltonian.

In recent years, this model has been extensively studied
and by now a wealth of results have been obtained~\cite{HPH,MR3922775,Bennett2017AddingRE,MR4025392,MR4130332,MR4052848,HPER,MR4052851,MR3595872}.
In this paper, we are interested in relaxing the (vanishing) minimum degree conditions in the results by a density condition, especially for the spanning subgraph problems.
Let us first give some known results.

Given $k$-graphs $F$ and $G$, an \emph{$F$-factor} of $G$ is a collection of vertex-disjoint copies of $F$ in $G$ that covers all the vertices of $G$.
The thresholds for $F$-factors in random graphs have been a difficult problem and are essentially known now.
Given a $k$-graph $F$, let $v_F$ and $e_F$ denote the number of vertices and edges of $F$, respectively.
If $F$ has at least two vertices, define
\[
m_1(F):=\max\left\{ d(F'):F'\subseteq F, v_{F'}\geq 2 \right\}, \text{ where } d(F'):=\frac{e_{F'}}{v_{F'}-1}.
\]
A breakthrough result on the threshold for $F$-factors was made by Johansson, Kahn, and Vu~\cite{factorth}, and now the thresholds are known to lie between $n^{-1/m_1(F)}$ and $n^{-1/m_1(F)}\log n$, and are known exactly for certain $F$.
For the $F$-factor problem in randomly perturbed graphs, Balogh, Treglown, and Wagner~\cite{MR3922775} proved the following result.

\begin{theorem}[\cite{MR3922775}]
\label{thm:BTW}
Let $F$ be a fixed graph with $e_F > 0$, and let $n \in v_F \mathbb{N}$.
For every $\alpha > 0$, there exists a constant $c = c(\alpha, F) > 0$ such that if $p \ge cn^{-1/m_1(F)}$ and $G_{\alpha}$ is an $n$-vertex graph with minimum degree at least $\alpha n$, then w.h.p.\ $G_{\alpha} \cup G(n, p)$ contains an $F$-factor.
\end{theorem}

Graphs with bounded maximum degree is another natural graph class to study.
More precisely, for any $n,\Delta\in \mathbb{N}$, let $\mathcal{F}(n,\Delta)$ be the collection of all $n$-vertex graphs with maximum degree at most $\Delta$.
After a series of works, Frankston, Kahn, Narayanan, and Park~\cite{MR4298747} showed that for every graph $G \in \mathcal{F}(n, \Delta)$, the threshold for $G$ is at most $n^{-2/(\Delta+1)} (\log n)^{2/(\Delta(\Delta+1))}$ (indeed achieved by $K_{\Delta+1}$-factors).

Addressing a problem in~\cite{MR3595872}, B\"ottcher, Montgomery, Parczyk, and Person~\cite{MR4130332} showed the following result for bounded-degree subgraphs in $G_{\alpha} \cup G(n,p)$.
\begin{theorem}[\rm \cite{MR4130332}]\label{boundeddegree}
Let $\alpha > 0$ be a constant and $\Delta \geq 5$ be an integer.
Suppose $G_{\alpha}$ is an $n$-vertex graph with minimum degree at least $\alpha n$.
Then for every $F \in \mathcal{F}(n, \Delta)$ and $p = \omega(n^{-2/(\Delta+1)})$, w.h.p.\ $G_{\alpha} \cup G(n,p)$ contains a copy of $F$.
\end{theorem}

Another class of extensively studied objects are uniform cycles in $k$-graphs, for which we need the following definitions.
For $k > \ell \geq 1$, we say a $k$-graph $H$ is a $k$-uniform \emph{$\ell$-cycle} if its vertices can be ordered cyclically $x_1, \ldots, x_n$ such that each edge contains $k$ consecutive vertices and any two subsequent edges intersect at exactly $\ell$ vertices.
Equivalently, $\{x_{(k-\ell)i+1}, \ldots, x_{(k-\ell)i+k}\}$ is an edge for every $i \in [n/(k-\ell)]$.
(When $\ell=k-1$, it is usually called a \emph{tight cycle}.)
Note that we need $(k-\ell) \mid n$ to ensure its existence.
For a hypergraph $H$, a \emph{Hamilton $\ell$-cycle} is a spanning subgraph of $H$ that is isomorphic to an $\ell$-cycle.
The thresholds for $\ell$-Hamiltonicity have been studied by Dudek and Frieze \cite{WOS:000287601500007,dudek2013tight},
 who proved that the asymptotic threshold is $1/n^{k-\ell}$ for $\ell\geq2$ and $\log n/n^{k-1}$ for $\ell=1$ (they also gave a
 sharp threshold for $k\geq4$ and $\ell=k-1$).
 It is also natural to consider $\ell$-Hamiltonicity in randomly perturbed $k$-graphs.
 In fact, Krivelevich, Kwan
 and Sudakov \cite{krivelevich2016cycles} extended the result of Bohman, Frieze and Martin \cite{MR1943857} to hypergraphs with minimum $(k-1)$-degree condition.
  Krivelevich, Kwan and Sudakov \cite{krivelevich2016cycles} asked whether the result in \cite{krivelevich2016cycles} can be extended to $\ell$-Hamiltonicity under
 minimum $d$-degree conditions for $d,\ell\in[k-1]$. McDowell and Mycroft \cite{WOS:000450299200005} found such results for
 $d\geq\max\{\ell,k-\ell\}$ and reiterated the question for arbitrary $d$ and $\ell$.
Finally, Han and Zhao \cite{han2020hamiltonicity} completely resolved this problem.
\begin{theorem}[\cite{han2020hamiltonicity}]
\label{thm:hamiltonicity}
    For integers $k\geq3,\ell\in[k-1]$ and $\alpha>0$ there exist $\varepsilon>0$ and an integer $c>0$
 such that the following holds for sufficiently large $n\in(k-\ell)\mathbb{N}$. Suppose $H$ is a $k$-graph on $n$ vertices with $\delta_1(H)\geq\alpha n^{k-1}$ and
 \begin{equation}
\nonumber
p=p(n)=
\left\{
\begin{aligned}
n^{-(k-\ell)-\varepsilon},\ {\rm if\ } \ell&\geq2, \\
cn^{-(k-1)},\ {\rm if\ }\ell&=1.
\end{aligned}
 \right.
\end{equation}
Then $H\cup G^{(k)}(n,p)$ contains a Hamilton $\ell$-cycle.
\end{theorem}

Given $k\geq2$ and $r\geq2$, we say that a $k$-graph is an $r$-\emph{th power of a tight cycle} if its vertices can be ordered cyclically such that each consecutive $k+r-1$ vertices span a copy of $K_{k+r-1}^{(k)}$ and there are no other edges than the ones forced by this condition.
For $k=2$ and $r\geq2$, the threshold for $r$-th power of a Hamilton cycle forms an interesting question, and is now known as $n^{-1/r}$, by Riordan~\cite{MR1762785} (for $r \ge 3$) and by Kahn, Narayanan and Park~\cite{MR4273128} (for $r=2$).
For $k\geq3$ and $r\geq2$, Parczyk and Person~\cite{parczyk2016spanning} proved that the threshold for the existence of the $r$-th power of a tight Hamilton cycle in $G^{k}(n,p)$ is $n^{-\binom{k+r-2}{k-1}^{-1}}$.
B\"ottcher, Montgomery, Parczyk and Person~\cite{MR4130332} studied power of Hamilton cycles in randomly perturbed graphs.
For $r$-th power of tight Hamilton cycles in $k$-graphs with $k\geq3$ and $r\geq2$, this was studied by Bedenknecht, Han, Kohayakawa and Mota~\cite{MR4025389} in randomly perturbed $k$-graphs and improved by Chang, Han and Thoma~\cite{chang2023powers}.
Let us state these results below.

\begin{theorem}
[\cite{MR4130332,chang2023powers}]
\label{thm:power_hamilton_cycle}
For $k\geq2, r \geq 2$ and constant $\alpha > 0$, there exists a constant $\varepsilon> 0$ such that if $G_{\alpha}$ is an $n$-vertex $k$-graph with minimum $\delta_{k-1}(G_\alpha)\geq\alpha n$ and $p= n^{-\binom{k+r-2}{k-1}^{-1} - \varepsilon}$, then w.h.p.\
$G_{\alpha} \cup G^{(k)}\bigl(n,p\bigr)$
contains the $r$-th power of a tight Hamilton cycle.
\end{theorem}

The assumptions in Theorems~\ref{thm:BTW}--\ref{thm:power_hamilton_cycle} are essentially best possible in the following sense.
Firstly, the minimum degree conditions on $G_\alpha$ are vanishing and clearly cannot be completely omitted;
secondly, as there are instances of $G_{\alpha}$ with an independent set of size $(1-\alpha)n$ (e.g. one may take $G_\alpha$ as a complete bipartite graph $K_{\alpha n, (1-\alpha)n}$), one has to embed most part of the spanning subgraph in $G^{(k)}(n,p)$, yielding that the assumptions on $p$ are essentially best possible (e.g., up to the value of $\varepsilon$ or $c$).

\subsection{Main results}
In this paper, we investigate relaxations of the above results, that is, replacing the minimum degree condition in the randomly perturbed $k$-graph model by a density condition.
Note that in some situations such a weakening is not possible.
For instance, for Hamiltonicity, in the result of Bohman, Frieze and Martin, one cannot replace the minimum degree condition by a density condition.
Indeed, if $G$ contains $\Omega(n)$ isolated vertices, then $G \cup G(n, C/n)$ will w.h.p. contain isolated vertices, and thus does not contain a Hamilton cycle.
Similarly, such a result is also impossible for spanning trees and for $F$-factors when $F$ is a forest without isolated vertices.

Our main results are strengthenings of Theorems~\ref{thm:BTW}--\ref{thm:power_hamilton_cycle}, showing that the minimum degree conditions in them can be replaced by a density condition, with two exceptions when $F$ is a forest in Theorem~\ref{thm:BTW} and when $\ell=1$ in Theorem~\ref{thm:hamiltonicity}.

\begin{theorem}\label{mainfactor}
    Let $F$ be a graph with $e_F > 0$ such that $F$ is not a forest, and let $n \in \mathbb{N}$ be divisible by $v_F$.
    For every $\varepsilon > 0$, there exists a constant $c = c(F, \varepsilon) > 0$ such that if $p \ge c n^{-1/m_1(F)}$ and $G$ is an $n$-vertex graph with at least $\varepsilon n^2$ edges, then w.h.p.\ $G \cup G(n, p)$ contains an $F$-factor.
\end{theorem}

\begin{theorem}
\label{mainbounded}
For every $\varepsilon > 0$ and $\Delta \geq 2$, there exists a constant $c = c(\varepsilon, \Delta)$ such that the following holds.
Let $F$ be an $n$-vertex graph with maximum degree $\Delta(F) \leq \Delta$.
If $p \geq c n^{-\frac{2}{\Delta + 1}}$ and $G$ is an $n$-vertex graph with at least $\varepsilon n^2$ edges, then w.h.p.\ $G \cup G(n, p)$ contains a copy of $F$.
\end{theorem}

We strengthen the result of Han and Zhao \cite{han2020hamiltonicity} for $\ell\geq2$ as follows.

\begin{theorem}
\label{hypergraph}
For integers $k\geq3$, $\ell\in[2,k-1]$ and $\varepsilon>0$, there exist $\eta>0$ and integer $C>0$ such that the following holds for sufficiently large $n\in(k-\ell)\mathbb{N}$.
Suppose $H$ is a $k$-graph on $n$ vertices with at least $\varepsilon n^k$ edges and $p=p(n)\geq n^{-(k-\ell)-\eta}$.
Then w.h.p., $H\cup G^{(k)}(n,p)$ contains a Hamilton $\ell$-cycle.
\end{theorem}

\begin{theorem}\label{mainpower}
For $k\geq2, r\geq2$ and $\varepsilon>0$, there exists $\eta>0$ such that the following holds.
Suppose $H$ is an $n$-vertex $k$-graph with at least $\varepsilon n^k$ edges and $p=p(n)\geq n^{-\binom{k+r-2}{k-1}^{-1}-\eta}$.
Then w.h.p. $H\cup G^{(k)}(n,p)$ contains the $r$-th power of a tight Hamilton cycle.
\end{theorem}

For Theorem~\ref{mainfactor}, we employ the ``absorption technique'' pioneered by R\"{o}dl, Ruci\'{n}ski and Szemer\'{e}di~\cite{MR2195584}, combined with several tools from the theory of binomial random graphs.
The key idea is to construct absorbers for any set of $v_F$ vertices.
To obtain an $F$-factor in $G \cup G(n,p)$, we follow a two-step strategy: first, we use the edges of $G$ to partially embed copies of $F$, and then we use the random edges in $G(n,p)$ to fill in the missing edges to yield a copy of $F$.
Theorem~\ref{hypergraph} and Theorem~\ref{mainpower} are proved by an application of the spread method, which has been successful in several recent breakthroughs.
Finally, we derive Theorem~\ref{mainbounded} by combining Theorem~\ref{mainfactor} with some known threshold results in $G(n,p)$.


\section{Notation and Preliminaries}

\subsection{Basic notations}
Given a graph $G$, and (not necessarily disjoint) subsets $S, S' \subseteq V(G)$, we write $E_G[S, S']$ for the set of edges in $G$ with one endpoint in $S$ and the other in $S'$.
We simply write $E[S, S']$ when the graph $G$ is clear from context.
We denote by $G[S]$ the subgraph of $G$ induced on the vertex set~$S$.
For any edge $e \in E(G)$, we write $G \cup e$ and $G \setminus e$ for the graphs obtained from $G$ by adding or deleting the edge~$e$, respectively.
For any vertex $v \in V(G)$, we denote by $G - v$ the graph obtained from $G$ by deleting $v$ along with all edges incident to it.

We use $\ll$ to denote a hierarchy between constants.
If we write that a statement holds whenever $0 < a \ll b, c \ll d$, it means that there exist non-decreasing functions $g_1, g_2 \colon (0,1] \to (0,1]$ and $f \colon (0,1]^2 \to (0,1]$ such that the statement holds for all $a, b, c, d$ satisfying $b \le g_1(d)$, $c \le g_2(d)$, and $a \le f(b, c)$.
We will not explicitly compute these functions to avoid cluttering the presentation of the proofs.

\subsection{Small subgraphs in a random graph}

We first introduce some notation from~\cite{MR1782847}.
For a graph $F$, let
\[
\Phi_F = \Phi_F(n, p) := \min\left\{ n^{v_{F'}} p^{e_{F'}} : F' \subseteq F,\ e_{F'} > 0 \right\}.
\]
For a graph $F$ and an independent set $W \subseteq V(F)$, define
\[
\Phi_{F, W} = \Phi_{F, W}(n, p) := \min\left\{ n^{v_{F'} - v_{F'[W]}} p^{e_{F'}} : F' \subseteq F,\ e_{F'} > 0 \right\}.
\]
It is helpful to think of $\Phi_F$ as the expected number of copies of $F$ in $G^{(k)}(n,p)$, and $\Phi_{F, W}$ as the expected number of rooted-copies of $F$, that is, with vertices in $W$ mapped to a fixed set of $|W|$ vertices.

It is easy to see that $\Phi_F = \Phi_{F, \emptyset}$ and $\Phi_{F \setminus W} \ge \Phi_{F, W}$ for any graph $F$ and any independent set $W \subseteq V(F)$.
The following result provides a sufficient condition that ensures $\Phi_F \ge cn$.
\begin{lemma}[{\cite[Lemma~2.4]{chang2022factors}}]
\label{otherfaiF}
    Let $c \ge 1$ be a constant, and let $F$ be a labelled graph.
    If $p = p(n) \ge c n^{-1/m_1(F)}$, then $\Phi_F \ge cn$.
\end{lemma}

One may expect that if $\Phi_F(n, p)$ is large, then the Erd\H{o}s-R\'enyi random graph $G(n, p)$ is very likely to contain a copy of $F$.
The following lemma confirms that this intuition is correct.
\begin{lemma}[{\cite[Lemma~2.2(i)]{MR4025389}}]
\label{otherFcover}
    Let $F$ be a labelled graph with $b$ vertices and $f$ edges. Suppose $1/n \ll 1/c \ll \lambda, 1/b, 1/f$.
    Let $G = G(n, p)$ be the binomial random graph.
    If $p = p(n)$ satisfies $\Phi_F(n, p) \ge cn$, then with probability at least $1 - \exp(-n)$, every induced subgraph of $G$ of order $\lambda n$ contains a copy of $F$.
\end{lemma}

Via a careful use of Janson's inequality, Han, Morris and Treglown~\cite{HPER}, further showed that if $\Phi_{F, W}$ is also large, then one can find an embedding of $F$ such that $W$ is mapped to a designated target set of vertices.
In this paper we use a version of this result stated in~\cite{chang2023powers}.

\begin{proposition}[{\cite[Proposition~5.2]{chang2023powers}}]\label{propembed}
 Let $k,b,f,w,\ell, n,t=t(n)\in \mathbb{N}$ be integers and suppose
 $1/n\ll 1/C, \gamma \ll \beta , 1/f, 1/b, 1/w, 1/\ell$.
Let $F$ be a labelled graph with base vertex set $W\subseteq V(F)$ which is an independent set of $F$ such that $|W|=w$, $|V(F)\setminus W|=b$, and $e_F=f$. Further, suppose that $p=p(n)$ satisfies $\Phi_{F\setminus W}\ge Cn$ and $\Phi_{F,W}\ge Cn^{1/\ell}$.

Let $V$ be an $n$-vertex set, and let $U_1,\ldots,U_t \subseteq V$ with $t\le \gamma n $
 be pairwise disjoint and be such that
$|U_i|=|W|$ for each $i\in [t]$. Suppose that $\mathcal{F}_1, \ldots, \mathcal{F}_t$ are families of ordered $b$-sets on $V$ such that $|\mathcal{F}_i|\ge \beta n^b$. Then w.h.p. there exists a collection of  embeddings $\phi_1, \ldots, \phi_t$ such that each $\phi_i$ embeds a copy of $F$ into $G^{(k)}(n,p)$ on $V$ with $W$  being mapped to $U_i$ and $V(F)\setminus W$ being mapped to a set in $\mathcal{F}_i$ which is vertex-disjoint with $\bigcup_{i\in [t]}U_i$. Furthermore, for $i\ne j$ we have $\phi_i(V(F)\setminus W)\cap \phi_j(V(F)\setminus W)=\emptyset$.
\end{proposition}

\subsection{Spreadness}
Recently, Frankston, Kahn, Narayanan, and Park~\cite{MR4298747} proved a fractional version of the Kahn-Kalai conjecture, proposed by Talagrand~\cite{Talagrand2010}.
Later, Park and Pham~\cite{MR4654612} established a stronger version of the Kahn--Kalai conjecture~\cite{MR2312440}.
An important corollary of the result by Frankston, Kahn, Narayanan, and Park is the connection between spreadness and thresholds in random graphs.

\begin{definition}[Spread]
Let $q \in [0, 1]$.
Assume that $\mathcal{H}$ is a hypergraph on the vertex set $V$, and let $\mu$ be a probability measure on the edge set of~$\mathcal{H}$.
We say that $\mu$ is \emph{$q$-spread} if for every subset $S \subseteq V$, the following holds:
\[
\mu\left(\left\{ A \in E(\mathcal{H}) : S \subseteq A \right\} \right) \leq q^{|S|}.
\]
\end{definition}

We use $V_p$ to denote a random subset of $V$ where each element $x \in V$ is included independently with probability~$p$.
We say that a hypergraph $\mathcal{H}$ is \emph{$r$-bounded} if every edge of $\mathcal{H}$ contains at most $r$ vertices.
\begin{proposition}[{\cite[Theorem~1.6]{MR4298747}}]
\label{FKNP}
Let $\mathcal{H}$ be an $r$-bounded hypergraph on vertex set $V$ that supports a $q$-spread distribution.
If $p \geq Kq \log r$, then with probability $1-o_r(1)$, the random subset $V_p$ contains an edge of $\mathcal{H}$.
\end{proposition}

In typical applications, the vertex set of the hypergraph $\mathcal{H}$ is the edge set of a graph $G$, and the hyperedges of $\mathcal{H}$ correspond to subgraphs of $G$ isomorphic to a graph~$H$.
Then $V_p$ represents the random set of edges of $G$, where each edge is included independently with probability~$p$.

Pham, Sah, Sawhney and Simkin~\cite{pham2023toolkit} introduced a notion of vertex-spreadness, which was further generalized by Kelly, M\"{u}yesser and Pokrovskiy~\cite{KELLY2024507}.
\begin{definition}[Vertex-spread]
Let $X$ and $Y$ be finite sets and let $\mu$ be a probability distribution over injections $\varphi: X \rightarrow Y$.
For $q \in [0, 1]$, we say that $\mu$ is a \emph{$q$-vertex-spread} if for every $s \le |X|$ and every pair of sequences of distinct vertices $x_1, \ldots, x_s \in X$ and $y_1, \ldots, y_s \in Y$, we have
\[
\mu\left( \left\{ \varphi : \varphi(x_i) = y_i \text{ for all } i \in [s] \right\} \right) \le q^s.
\]
\end{definition}

A \emph{hypergraph embedding} $\psi : G \rightarrow H$ of a hypergraph $G$ into a hypergraph $H$ is an injective map $\psi : V(G) \rightarrow V(H)$ that maps edges of $G$ to edges of $H$.
Thus, there exists an embedding of $G$ into $H$ if and only if $H$ contains a subgraph isomorphic to $G$.
Note that a uniformly random embedding $\psi : G \rightarrow H$, when $G$ and $H$ have the same vertex set and $H$ is a complete hypergraph, corresponds to a random permutation of $V(H)$.
Such an embedding is $(e/|V(H)|)$-vertex-spread (by Stirling's approximation).
The following result allows us to connect spreadness and vertex-spreadness.
\begin{proposition}[{\cite[Proposition~1.17]{KELLY2024507}}]
\label{kelly}
For every $k, \Delta \in \mathbb{N}$ and $C > 0$, there exists a constant $C_{\ref{kelly}} = C_{\ref{kelly}}(C, k, \Delta) > 0$ such that the following holds for all sufficiently large $n$.

Let $H$ and $G$ be $n$-vertex $k$-graphs with $\Delta(G) \leq \Delta$.
If there exists a $(C/n)$-vertex-spread distribution on embeddings $\psi: G \rightarrow H$, then there exists a $(C_{\ref{kelly}} / n^{1/m_1(G)})$-spread distribution on subgraphs of $H$ isomorphic to $G$.
\end{proposition}

We note that Proposition~\ref{kelly} implies that if $G$ is a $K_{d+1}$-factor, then there exists a $(C / n^{2/(d+1)})$-spread distribution on subgraphs of $K_n$ that are isomorphic to $G$.
Frankston, Kahn, Narayanan and Park~\cite{MR4298747} further showed that if $G$ is far from being a $K_{d+1}$-factor, then one can obtain a significantly better spread distribution. Combining their result with Proposition~\ref{FKNP}, we obtain the following.
\begin{lemma}[{\cite[Lemma~7.6]{MR4298747}}]
\label{FKNP7.6}
There exists a constant $\varepsilon_{\ref{FKNP7.6}} = \varepsilon_{\ref{FKNP7.6}}(d) > 0$ such that the following holds.
Let $H$ be an $n$-vertex graph with $\Delta(H) \leq d$, and suppose $H$ has no connected component isomorphic to $K_{d+1}$.
If $p \ge n^{-2/(d+1) - \varepsilon_{\ref{FKNP7.6}}}$, then $G(n, p)$ contains a copy of $H$ w.h.p.
\end{lemma}

\section{$F$-factors in randomly perturbed graphs}
\subsection{Outline of Theorem~\ref{mainfactor}}
Let $F$ be a graph with $m$ vertices that is not a forest.
To find an $F$-factor in $G \cup G(n, p)$, we employ the absorption method.
In most applications of the absorption method, one finds a set called an absorbing set that can `absorb' any small set of vertices.
More precisely, one finds a set $V_0 \subseteq V(G)$ such that for any leftover vertices $R \subseteq V(G) \setminus V_0$ of sufficiently small size, $G[V_0 \cup R]$ contains an $F$-factor.
Thus, by showing that $G \setminus V_0$ contains an almost $F$-factor, one can conclude the proof.
On the other hand, in our case, we cannot find a set that can absorb every small set of vertices.
Instead, we can find a set $V_0$ such that for each small set $R$, $V_0$ can absorb $R$ with high probability.
To discuss the existence of such an absorbing set, we first define an `absorber'.
\begin{definition}
    Let $F$ be a graph with $m$ vertices, and let $H$ be a graph with $n$ vertices.
    For a set $S \subseteq V(H)$ with $|S|=m$ and a set $A \subseteq V(H) \setminus S$, we say that $A$ is an \emph{$S$-absorber (for $F$)} if both $H[A]$ and $H[A \cup S]$ contain an $F$-factor.
\end{definition}

In order to obtain an absorbing set, we need the following result for building the absorbing set.
It builds on the bipartite-template method of Montgomery~\cite{MR3998769} and several versions of such a result already exist in the literature (see e.g.~\cite{chang2022factors, Nenadov2020}).
Roughly speaking, in these versions it is assumed that (w.h.p.)~every $m$-vertex set has linearly many vertex-disjoint absorbers.
However, for our problem we could only prove weaker assumptions so we provide here an appropriate version for our use.
Indeed, in comparison, what we can show is that given a set of linearly many vertex-disjoint $m$-vertex sets, w.h.p.~one can find vertex-disjoint absorbers one for each of them.

We include a proof of this lemma in the appendix for completeness and remark that in the statement below, by saying that $G^*$ is a \textit{random graph}, we just mean that $G^*$ is drawn from certain probability distribution (so not necessarily binomial random graphs), which allows the lemma to be applicable under other contexts.

\begin{lemma}\label{changabsorbplus}
    Let $F$ be a  graph with $m$ vertices. Suppose $1/n\ll \xi \ll \gamma_1\ll \tau \ll \gamma \ll 1/m$ and let $C_1>0$ be a constant.
   For any $n$-vertex random graph $G^*$, if there is a constant $C_1>0$ such that $G^*$ satisfies the following  properties:
   \begin{enumerate}[$(i)$]
    \item  There is a vertex set $V_1\subseteq V(G^*)$ with  $\tau n/4 \le|V_1|\le \tau n$. \label{cond:absorbing1}

 	\item
  For any $S\subseteq V(G^*)$ with $|S|\le \gamma_1 n$,  w.h.p.  there exists a set $T\subseteq V_1\setminus S$ with $|T|=(m-1)|S|$ such that $G^*[S\cup T]$ contains an $F$-factor. \label{cond:absorbing2}
	
 \item Let  $k\le \frac{120\tau}{\gamma_1}$ be a constant and  $\mathcal{S}_1, \ldots,\mathcal{S}_k $ be  families of $m$-sets of  $V(G^*)$ such that each $\mathcal{S}_i$ is pairwise disjoint.
    Let $\mathcal{S}_i=\{S_{i,1},\ldots,S_{i,t_i}\}$, where $t_i=|\mathcal{S}_i|\le \gamma_1 n$. Then w.h.p.~there exist vertex-disjoint families of $C_1$-sets $\{T_{i,j}|  i\in [k] , j\in [t_i]\}$ in $V(G^*)\setminus \bigcup_{i\in [k], j\in [t_i]} S_{i,j}$
 such that for any $i\in [k]$ and $j\in [t_i]$, $ T_{i,j}$ is an $S_{i,j}$-absorber. \label{cond:absorbing3}
\end{enumerate}

   Then $G^*$ contains an absorbing set $V_0\subseteq V(G^*)$ of size at most $\gamma n$ such that, for every $R\subseteq V(G^*)\setminus V_0$ with $|R|\le \xi n$ such that $m$ divides $|V_0|+|R|$, w.h.p. $G^*[V_0\cup R]$ contains an $F$-factor.
\end{lemma}

We use the following lemma to obtain a sufficient condition for the existence of an absorbing set.

\begin{lemma}\label{absorblem}
 	Let $F$ be a graph with $m$ vertices that is not a forest and $C_1=\lceil m_1(F) m^2\rceil m^2$. Let $1/n\ll \gamma_1  \ll \tau \ll \varepsilon, 1/m$.
    Then  there is a constant $c>0$ such that
    if $p\ge cn^{-\frac{1}{m_1(F)}}$ and $e(G)\ge \varepsilon n^2$, then $G^*:=G\cup G(n,p)$  contains a vertex set $V_1\subseteq V(G^*)$ satisfying the properties \ref{cond:absorbing1}-\ref{cond:absorbing3} in Lemma~\ref{changabsorbplus}.
\end{lemma}

Applying Lemmas~\ref{changabsorbplus} and \ref{absorblem}, we  prove Theorem~\ref{mainfactor} as follows.

\begin{proof}[Proof of Theorem~\ref{mainfactor}]
Suppose $F $ is a graph with $m$ vertices and $F$ is not a forest. Let $n\in m\mathbb{N}$. Let $\varepsilon >0$ be given and suppose $1/n \ll \xi \ll \gamma \ll \varepsilon,1/m$.
Let $c'$ play the role of  $c$ in Lemma \ref{absorblem} and  let $0<1/c\ll 1/c'$.

Suppose $p\ge cn^{-\frac{1}{m_1(F)}}$ and let $G$ be an $n$-vertex graph with vertex set $V$ such that $e_G\ge \varepsilon n^2$. We will expose $G':=G(n,p)$ in two rounds: $G'=G_1\cup G_2$, where $G_1$ and $G_2$ are independent copies of $G(n,p')$ and $(1-p')^2=1-p$.  Since $(1-p')^2>1-2p'$, $p'>\frac{p}{2}\ge c'n^{-\frac{1}{m_1(F)}}$.

Let $G^1=G\cup G_1$ for simplicity.
By Lemmas \ref{changabsorbplus} and \ref{absorblem}, we obtain that $G^1$ contains an absorbing set $V_0$ of size at most $\gamma n$ such that for every $R\subseteq V\setminus V_0$ with $|R|\le \xi n$ such that $m$ divides $|V_0|+|R|$,  w.h.p. $G^1[V_0\cup R]$ contains an $F$-factor.

Next, we aim to find an almost $F$-factor in $G_2$ that covers all of $V\setminus V_0$ but at most $\xi n$ vertices. As $p'>c'n^{-\frac{1}{m_1(F)}}$, Lemma \ref{otherfaiF} gives   $\Phi_F(n,p')\ge c'n$.
Thus, by Lemma \ref{otherFcover} with $\lambda =\xi$, w.h.p. every induced subgraph of $G_2$ of order $\xi n$ contains a copy of $F$.
By selecting vertex-disjoint copies of $F$ greedily, we find pairwise vertex-disjoint copies of $F$ in $V\setminus V_0$ that covers all but at most $\xi n$ vertices.
Denote by $R$ the set of the remaining vertices of $V\setminus V_1$. Then as $|R|\le \xi n$, we conclude that w.h.p. $G^1[V_0\cup R]$ contains an $F$-factor.
Therefore, w.h.p. $G\cup G_1\cup G_2$ contains an $F$-factor.
\end{proof}

To prove Lemma~\ref{absorblem}, we need to find an $S$-absorber robustly for a set $S$ of size $m$ in $G \cup G(n, p)$.
The rest of this section is devoted to proving Lemma~\ref{absorblem}.

\subsection{Absorbing}
To prove Lemma~\ref{absorblem}, we construct an absorbing structure by partitioning the edges into two groups, so that one group will be found in $G$ and the other part will be found in $G(n,p)$.
To do this, we first prove the following proposition, which characterizes properties of the minimal subgraphs $F' \subseteq F$ satisfying $\frac{e_{F'}}{v_{F'} - 1} = m_1(F)$.
\begin{proposition}\label{minisubg}
    Given a graph $F$, and for $u \in V(F)$, let
    $\mathcal{F}'_u = \{ F' \subseteq F : u \in V(F'), \ \frac{e_{F'}}{v_{F'} - 1} = m_1(F) \}$
    and let $\mathcal{F}_u$ be the set of all minimal elements of $\mathcal{F}'_u$ with respect to inclusion. Then the following hold:
    \begin{enumerate}[$(i)$]
        \item For every pair of distinct subgraphs $F_1, F_2 \in \mathcal{F}_u$, we have $V(F_1) \cap V(F_2) = \{u\}$. \label{prop:minisubg1}
        \item If $d_F(u) = \delta(F)$, then $|\mathcal{F}_u| \leq 1$. \label{prop:minisubg2}
          \item Let $u$ be a vertex with $d_F(u) = \delta(F)$. If $|\mathcal{F}_u| = 0$, then there exists an edge in $F$ with at most one end in $N_F[u]$.  \label{prop:minisubg3}
          \item If $F$ is not a forest, then  any $F_u\in \mathcal{F}_u$ is not a tree. \label{prop:minisubg4}
    \end{enumerate}
\end{proposition}
\begin{proof}
    (i): By definition, it is clear that $u\in V(F_1)\cap V(F_2)$.
    Assume, for contradiction, that $|V(F_1) \cap V(F_2)| \geq 2$.
    Let $F_3 = F_1 \cap F_2$, where $V(F_3) = V(F_1) \cap V(F_2)$ and $E(F_3) = E(F_1) \cap E(F_2)$.
    By the minimality of $F_1$, we have  $\frac{e_{F_3}}{v_{F_3}-1}<m_1(F)$.
    Since $e_{F_i} = m_1(F) \cdot (v_{F_i} - 1)$ for $i \in [2]$ and $e_{F_3} < m_1(F) \cdot (v_{F_3} - 1)$, we  obtain that
    \begin{align*}
        m_1(F)&\ge \frac{e_{F_1\cup F_2}}{v_{F_1\cup F_2}-1}
        \\&=\frac{e_{F_1}+e_{F_2}-e_{F_3}}{v_{F_1}+v_{F_2}-v_{F_3}-1}\\
        &>\frac{m_1(F)\cdot\big[v_{F_1}-1+v_{F_2}-1-\big(v_{F_3}-1\big)\big]}{v_{F_1}+v_{F_2}-v_{F_3}-1}\\
        &=m_1(F),
    \end{align*}
    which is a contradiction. Thus, $V(F_1)\cap V(F_2)=\{u\}$.

    (ii): Assume that $\mathcal{F}_u$ contains at least two elements.
    Let $F_1, F_2 \in \mathcal{F}_u$ be two distinct subgraphs.
    By part (i), there exists $i\in [2]$ such that $d_{F_i}(u)\le \frac{\delta(F)}{2}$.
    If $v_{F_i} = 2$, then $F_i$ consists of a single edge, and so $m_1(F) \leq 1$.
    In this case, $F$ is a forest.
    Since $u$ is a vertex of minimum degree, we have $d_F(u) = 1$, and by (i), it follows that $|\mathcal{F}_u| \leq 1$, contradicting our assumption.
    So we may assume that $v_{F_i}\ge 3$.
    Observe that
    \[m_1(F)\ge m_1(F_i-u) = \frac{e_{F_i-u}}{v_{F_i-u}-1}=\frac{e_{F_i-u}}{v_{F_i}-2}.\]
    On the other hand, by the definition of $m_1(F_i) = m_1(F)$, we also have
    \begin{align*}
        m_1(F)&=\frac{e_{F_i}}{v_{F_i}-1}\\
        &=\frac{e_{F_i-u}+d_{F_i}(u)}{v_{F_i}-1}\\
        &\le  \frac{ m_1(F)\cdot \big( v_{F_i}-2\big)+d_{F_i}(u)}{v_{F_i}-1}\\
        &=m_1(F)+\frac{d_{F_i}(u)-m_1(F)}{v_{F_i}-1}.
    \end{align*}
    It yields that $\frac{\delta(F)}{2}-m_1(F)\ge d_{F_i}(u)-m_1(F)\ge 0$.
    However, from the definition of $m_1(F)$, we have \begin{align*}
    m_1(F)&\ge \frac{e_F}{v_{F}-1}\ge \frac{ \delta(F)\cdot v_{F}}{2\big(v_{F}-1\big)}>\frac{\delta(F)}{2},
    \end{align*}
    which contradicts the inequality above.
    Therefore, $|\mathcal{F}_u|\le 1$.

    (iii): Suppose that  every edge of $F$ has both ends in $N_F[u]$. Note that  $d_F(u) = \delta(F)$. Then  $F$ is a clique. But this implies that $|\mathcal{F}_u| = 1$, a contradiction.

    (iv): This follows immediately from the fact that $m_1(F) > 1$ for non-forest graphs.
\end{proof}

According to Proposition~\ref{minisubg}~(ii), we give the following definition, which will be useful for embedding $S$-absorbers in $G \cup G(n,p)$.

\begin{figure}[htbp]
	\centering \includegraphics[scale=1]{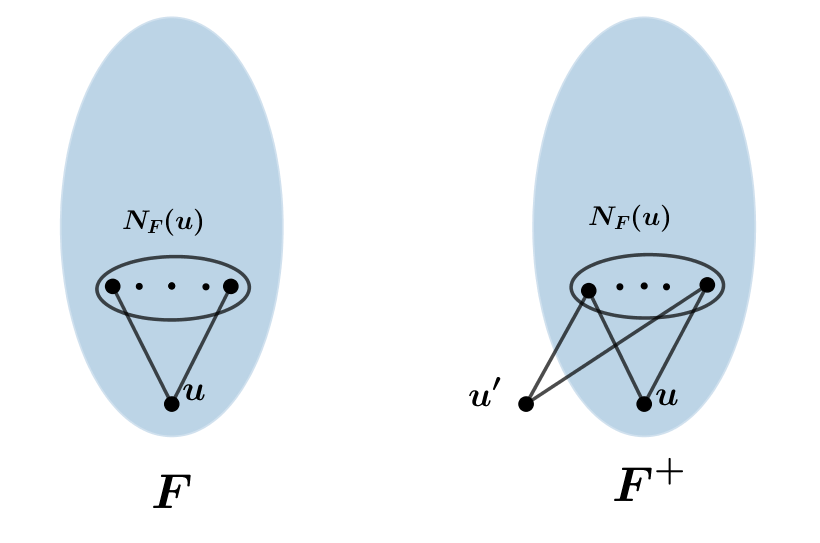}\vspace{-5mm}
	\caption{The graphs $F$ and $F^+$.}
	\label{FandFplus}
\end{figure}

\begin{figure}[h!]
    \begin{minipage}[b]{0.48\textwidth}
        \centering
        \includegraphics[width=0.8\textwidth]{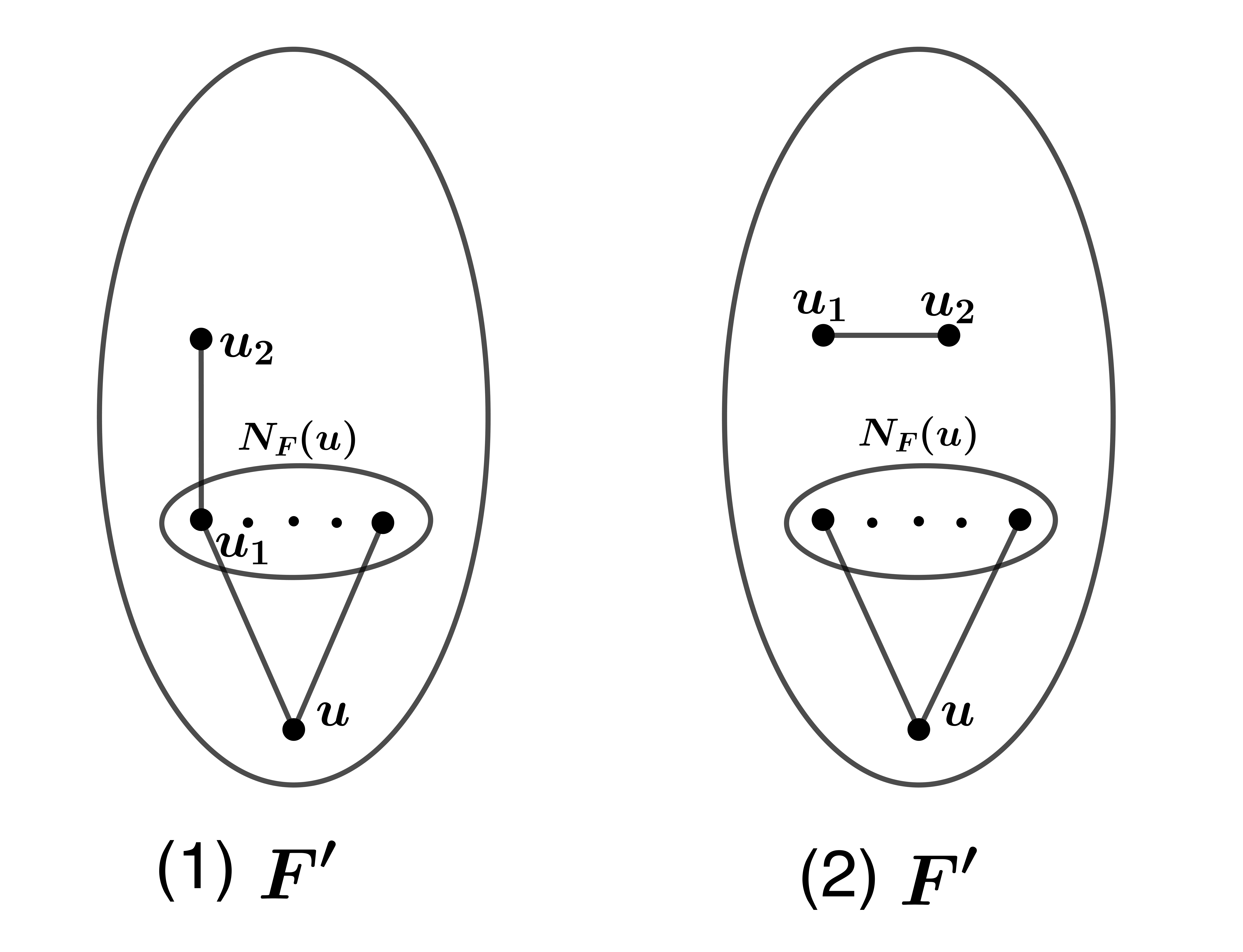}
        \caption{Two types of $F'$ in (i).}
        \label{Fi}
    \end{minipage}
    \hspace{0.5cm}
    \begin{minipage}[b]{0.45\textwidth}
        \centering
        \includegraphics[width=0.68\textwidth]{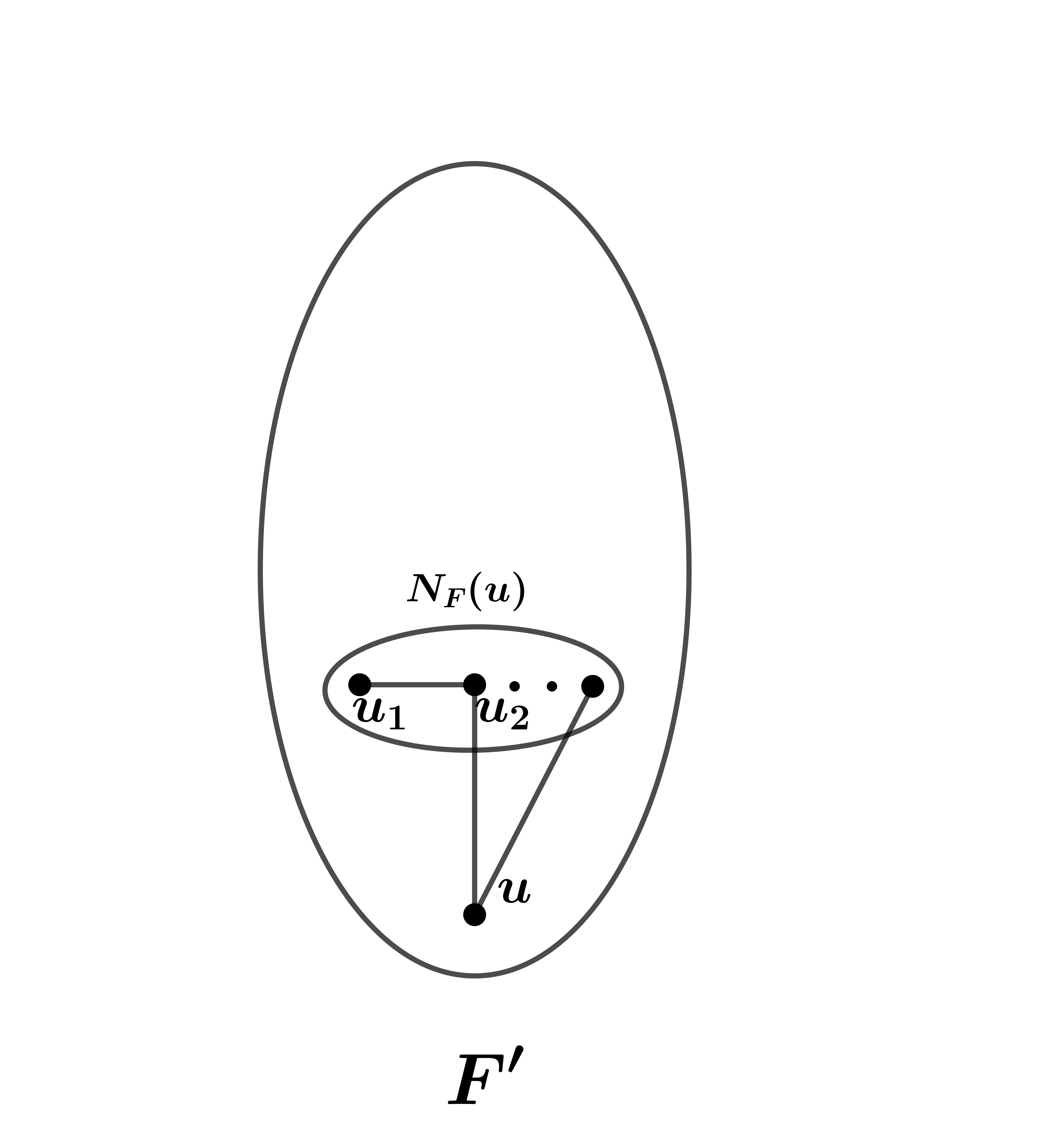}
        \caption{The graph $F'$ in (ii).}
         \label{Fii}
    \end{minipage}
\end{figure}

\begin{definition}\label{FinGnp} \rm{($u,u',F_u, F',F^*,u_1u_2$).}
    Let $F$ be a graph that is not a forest. Let $u \in V(F)$ be a vertex of minimum degree and let $F^+$ be the graph obtained by creating a new copy of $u$, denoted $u'$, as illustrated in Figure~\ref{FandFplus}. We now define $F'$ and $F^*$ as follows:

    \begin{itemize}
        \item[(i)] If $|\mathcal{F}_u|=0$, then by Proposition \ref{minisubg} (iii), there exists an edge with at most one end  in $N_F[u]$. In this case, there is an edge $u_1u_2 \in E(F)$ such that $u_i \notin N_F[u]$ for some $i \in [2]$.
        If $|\mathcal{F}_u|=1$, say $\mathcal{F}_u = \{F_u\}$, and
        $F_u$ contains an edge with at most one end in $N_{F_u}[u]$, then there is an edge $u_1u_2 \in E(F_u)$ such that $u_i \notin N_F[u]$ for some $i \in [2]$. \\
        Both of the cases, define $F'\subseteq F$ as the forest with $V(F')=V(F)$ and $E(F')=E[u,N_{F}(u)] \cup \{u_1u_2\}$, see Figure~\ref{Fi} for the illustration.
        Define $F^*$ as the graph obtained from $F^+$ by deleting all edges in $F'$, i.e.,
        $V(F^*)=V(F^+)$ and $E(F^*)=E(F^+)\setminus E(F')$.

    \item[(ii)] If $|\mathcal{F}_u|=1$, say $\mathcal{F}_u = \{F_u\}$, and every edge of  $F_u$
    has both ends in $N_{F_u}[u]$,
    then, since $F_u$ is not a forest, Proposition \ref{minisubg} (iv) implies that there exists an edge $u_1u_2 \in E(F_u)$ such that $u_1, u_2 \in N_{F_u}(u)$. \\
    Define $F'$ as the forest with $V(F') = V(F)$ and
    $E(F')=E[u,N_{F}(u)\setminus\{u_1\}] \cup \{u_1u_2\}$,see Figure~\ref{Fii} for the illustration.
    Define $F^*$ as the graph obtained from $F^+$ by deleting all edges in $F'$, i.e.,
    $V(F^*)=V(F^+)$ and $E(F^*)=E(F^+)\setminus E(F')$.
    \end{itemize}
\end{definition}

For the rest of this section, we fix a graph $F$ that is not a forest, and a vertex $u \in V(F)$ of minimum degree.
The following lemma provides lower bounds for $\Phi_{F^* \setminus W}$ and $\Phi_{F^*, W}$.

\begin{lemma}\label{Ffaivalue}
Let $F^*$ be the labelled graph defined in Definition~\ref{FinGnp},  $\varepsilon_0 = \frac{1}{m_1(F)v_F^2}$, and $W = \{u'\}$.
For every $C > 0$, there exists a constant $c > 0$ such that the following holds:
if $p \ge cn^{-\frac{1}{m_1(F)}}$, then
$\Phi_{F^*\setminus W}\ge  Cn$ and $\Phi_{F^*,W}\ge  Cn^{\varepsilon_0}$.
Furthermore, $m_1(F^*)\le m_1(F)$.
\end{lemma}
\begin{proof}
We choose $c > 0$ such that $0 < \frac{1}{c} \ll \frac{1}{C}$.
Let $p = cn^{-1/m_1(F)}$, and let $u_1u_2$ be the special edge of $F$ as in Definition~\ref{FinGnp}.
Let $H$ be a subgraph of $F^*$ with $v_H$ vertices and $e_H > 0$ edges.
By the definitions of $\Phi_{F^* \setminus W}$ and $\Phi_{F^*, W}$, it suffices to show that $n^{v_H}p^{e_H}\ge Cn$ if $u'\notin V(H)$ and $n^{v_H-1}p^{e_H}\ge Cn^{\varepsilon_0}$ if $u'\in V(H)$.
When $u'\notin V(H)$, we have $H\subseteq F$.
It follows that $\frac{e_H}{v_H-1}\le m_1(F)$ and so $n^{v_H}p^{e_H}\ge c^{e_H} n^{v_H} n^{-\frac{e_H}{m_1(F)}}\ge Cn$.
It remains to consider the case where $u'\in V(H)$.
If we show $\frac{e_H}{v_H-1}< m_1(F)$  for all such $H$, then
$m_1(F^*) \le  m_1(F)$.
Moreover, as $m_1(F)$ is achieved by a subgraph $F^\circ$ with at most $v_F$ vertices, we indeed have $m_1(F) - m_1(F^*) = \frac{e_{F^\circ}}{v_{F^\circ}-1} - \frac{e_H}{v_H-1} \ge \frac{1}{v_F^2}$.
As $\varepsilon_0= \frac{1}{m_1(F) v^2_F}$, we get
\[
\frac{e_H}{v_H-1}\le m_1(F)-\frac{1}{v^2_F}= m_1(F) (1-\varepsilon_0).
\]
For any $H\subseteq F^*$, if $u'\in V(H)$, then $\frac{e_H}{v_H-1} \le m_1(F) (1-\varepsilon_0)$ and so
\[
n^{v_H-1}p^{e_H}\ge c^{e_H} n^{v_H-1} n^{-\frac{e_H}{m_1(F)}}\ge c^{e_H} n^{\varepsilon_0(v_H-1)}\ge Cn^{\varepsilon_0}.
\]
Therefore, it remains to show $d(H)=\frac{e_H}{v_H-1}< m_1(F)$ for every $ H \subseteq F^*$ with $u'\in V(H)$.

\medskip
{\noindent \bf Case 1. $u\notin V(H)$.  }
\medskip

By Definition~\ref{FinGnp}, the graph $F^*-u$ is isomorphic to the graph obtained from $F$ by deleting the edge $u_1u_2$, i.e.,
$F^*-u \cong  F\setminus \{u_1u_2\}$. Moreover, under this isomorphism, the vertex $u'$ in $F^* - u$ corresponds to $u$ in $F \setminus \{u_1u_2\}$.
Hence, $H$ is isomorphic to a subgraph of $F\setminus \{u_1u_2\}$ that contains $u$.
If $|\mathcal{F}_u| = 0$, then clearly $\frac{e_H}{v_H - 1} < m_1(F)$.
Otherwise, $|\mathcal{F}_u| = 1$ and by Definition~\ref{FinGnp}, the edge $u_1u_2$ belongs to $F_u \in \mathcal{F}_u$. In this case, since $u_1u_2 \notin E(H)$, it follows that $H$ is not isomorphic to $F_u$, and thus $d(H)=\frac{e_H}{v_H-1} < m_1(F)$.

\medskip
{\noindent \bf Case 2. $u\in V(H)$. }
\medskip

Recall that $d_{F^*}(u) \le 1$, and since $F$ is not a forest, we have $m_1(F)>1$.
If $e_H\le v_H-1$, then clearly $\frac{e_{H}}{v_{H}-1}\le 1 <m_1(F)$.
If $e_H\ge v_H$, let $H_1=H-u$, which is isomorphic to a subgraph of $F$ as in case 1. Then $\frac{e_{H}}{v_{H}-1}< \frac{e_{H}-1}{v_{H}-2}\le \frac{e_{H_1}}{v_{H_1}-1}   \le m_1(F)$.

Therefore, $\frac{e_H}{v_H-1} < m_1(F)$ for every $H \subseteq F^*$ with $u' \in V(H)$, which concludes the proof.
\end{proof}

\begin{figure}[htp]
    \centering
   \includegraphics[width=0.5\linewidth]{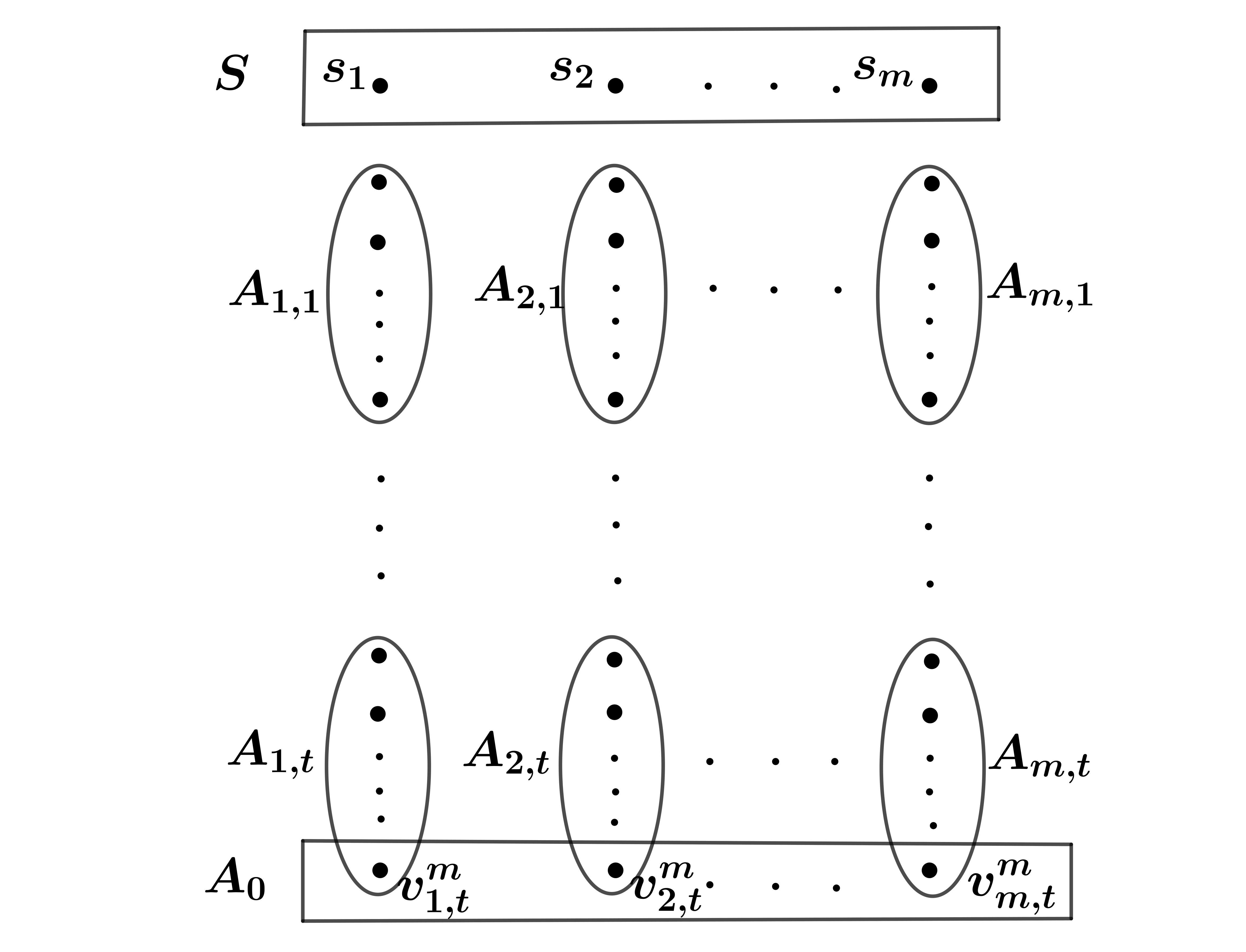}
    \caption{The vertex sets $A_{i,j}$, $S$ and $A_0$.}
    \label{someAij}
\end{figure}

We now define the absorbing structure that will be used in this section.
To find $S$-absorbers in $G \cup G(n,p)$, we divide the edges of absorbers into two parts: one part consists of edges we find in $G(n,p)$, and the other part consists of edges we find in $G$.
Formal definitions of these two parts are given in the following two definitions:
\begin{definition}\label{manyF}
Let $m = v_F$, and let $F^*$ be the graph defined in Definition~\ref{FinGnp}, with two special vertices $u$ and $u'$.
For each $i \in [m]$ and $j \in [t]$, define the set $A_{i,j} = \{v_{i,j}^1, \ldots, v_{i,j}^m\}$, where the sets $A_{i,j}$ are pairwise disjoint (see Figure~\ref{someAij}).
Let $A_0 = \{v_{1,t}^m, \ldots, v_{m,t}^m\}$, and let $S = \{s_1, \ldots, s_m\}$ be a set disjoint from all the vertex sets $A_{i,j}$.

Let $B$ be the graph with the vertex set $V(B)=\big(\bigcup_{i\in [m], j\in [t]} A_{i,j}\big) \cup S$ and the edge set is given by the following.
\begin{itemize}
    \item Let $B[A_0] \cong F$.
    \item For each $i \in [m]$, let $B[\{s_i\} \cup A_{i,1}] \cong F^*$, where $s_i$ and $v_{i,1}^m$ correspond to $u'$ and $u$, respectively.
    \item For each $i \in [m]$ and $j \in [t-1]$, let $B[\{v_{i,j}^m\} \cup A_{i,j+1}] \cong F^*$, where $v_{i,j}^m$ and $v_{i,j+1}^m$ correspond to $u'$ and $u$, respectively.
\end{itemize}
The edge set of $B$ is then the union of the above components, i.e.,
$B = B[A_0]\cup \big(\bigcup_{i\in [m]} B[\{s_i\}\cup A_{i,1}] \big) \cup \big(\bigcup_{i\in [m],j\in [t-1]} B[\{v_{i,j}^m\}\cup A_{i,j+1}] \big)$.
\end{definition}

 \begin{definition}\label{manyF'}
Let $m = v_F$, and let $F'$ be the graph defined in Definition~\ref{FinGnp}.
For each $i \in [m]$ and $j \in [t]$, define the set $A_{i,j} = \{v_{i,j}^1, \ldots, v_{i,j}^m\}$, where the sets $A_{i,j}$ are pairwise disjoint (see Figure~\ref{someAij}).
Let $D$ be the forest with vertex set $V(D)=\big(\bigcup_{i\in [m], j\in [t]} A_{i,j}\big) $ and the edge set is given by the following.
For each $i \in [m]$ and $j \in [t]$, define $D[A_{i,j}] \cong F'$, where the isomorphism follows the same correspondence as in Definition~\ref{manyF}.
In particular, under this isomorphism, the vertex $v_{i,j}^m$ corresponds to $u$.
The edge set of $D$ is then given by $\bigcup_{i\in [m],j\in [t]} D[ A_{i,j}] $.
\end{definition}

Then we have that the union of $B$ and $D$ forms an absorber for $F$.
\begin{lemma}\label{prop:BD_is_absorber}
    Let $B$ and $D$ be the labelled graphs defined in Definitions~\ref{manyF} and~\ref{manyF'}, respectively.
    We consider $B$ and $D$ as graphs on the same vertex set by identifying $A_{i,j}$ in both $B$ and $D$ for each $i \in [m]$, $j \in [t]$.
    Then both $B \cup D$ and $(B \cup D) \setminus S$ contain an $F$-factor.
\end{lemma}
\begin{proof}
    By definition, the union $F^* \cup F'$ is isomorphic to $F^+$.
    Thus, for each $i \in [m]$ and $j \in [t]$, we have $(B \cup D)[A_{i,j}] \cong F$, which implies that $(B \cup D) \setminus S$ has an $F$-factor.
    Similarly, as $F^+ \setminus u \cong F$, it follows that $(B \cup D)[A_{i, j+1} \cup \{u_{i, j}^m\} \setminus \{u_{i, j+1}^m\}] \cong F$ and $(B \cup D)[A_{i, 1} \cup \{s_i\}] \cong F$ for each $i\in [m]$, $j\in [t-1]$.
    In addition, $B[A_0] \cong F$.
    Therefore, we  construct an $F$-factor in $B\cup D$ by using the edges in $B[A_0]$, $B[\{s_i\}\cup A_{i,1}]$, and $B[\{v_{i,j}^m\}\cup A_{i,j+1}]$ for each $i\in [m]$, $j\in [t-1]$.
\end{proof}

We now prove that many copies of $B \cup D$ can be found, well-distributed, in $G \cup G(n,p)$.
To do this, we first show that the labelled graph $B$ can be embedded robustly in $G(n,p)$.
\begin{proposition}{\rm (\cite{chang2022factors}, Proposition 2.2)}.\label{chang2phi}
    Let $F_1$ and $F_2$ be labelled graphs with $V(F_1)\cap V(F_2)=\{v\}$. Then $\Phi_{F_1\cup F_2}\ge \min\{\Phi_{F_1},\Phi_{F_2}, \Phi_{F_1}\Phi_{F_2}n^{-1}\}$.
\end{proposition}

Using Proposition~\ref{chang2phi}, we  prove the following lemma, which provides lower bounds for $\Phi_{B \setminus S}$ and $\Phi_{B, S}$.

\begin{lemma}\label{bigFfaivalue}
Let $B$ be the labelled graph defined in Definition~\ref{manyF} and $ \varepsilon_0 = \frac{1}{m_1(F) v^2_F}$.
For every $C>0$, there exists $c>0$ such that the following holds.
If $p\ge cn^{-\frac{1}{m_1(F)}}$, then $\Phi_{B\setminus S}\ge Cn$ and $\Phi_{B,S}\ge Cn^{\varepsilon_0}$.
\end{lemma}
\begin{proof}
    Let $1/c \ll 1/C$ be constants and we may assume that $C>1$.
    For a subgraph $H$ of $B$, let $\lambda(H)=n^{v_{H}}p^{e_H}$ and $\lambda'(H)=n^{v_{H}-v_{H[S]}}p^{e_H}$.
    We first prove that $\Phi_{B\setminus S}\ge Cn$.
    By Lemma~\ref{Ffaivalue}, we have $m_1(F^*) \le m_1(F)$.
    By Lemma~\ref{otherfaiF}, we know that $\Phi_F \ge cn \ge Cn$ and $\Phi_{F^*} \ge cn \ge Cn$.
    Since the graph $B$ is constructed by starting from a copy of $F$ on $A_0$ and then iteratively attaching copies of $F^*$, we  apply Proposition~\ref{chang2phi} inductively to obtain $\Phi_{B\setminus S} \geq \Phi_B \geq Cn$.

    Next, we prove that $\Phi_{B, S} \ge Cn^{\varepsilon_0}$.
    Let $H \subseteq B$ be any subgraph.
    We proceed by analyzing cases based on the number of vertices in $H \cap S$.
    Note that when $H_1, \ldots, H_k$ are connected components of $H$, then $\lambda'(H) = \prod_{i=1}^k \lambda'(H_i)$ so it sufficient to prove that $\lambda'(H)\ge Cn^{\varepsilon_0}$ for a connected subgraph $H$.
    If $|V(H)\cap S|=0$, then $H$ is a subgraph of $B\setminus S$ and by the  earlier result, we have $\lambda(H) = \lambda'(H) \ge Cn$.
    Next assume $|V(H)\cap S|=1$ and let $s_i$ be the unique vertex in $S \cap V(H)$.
    Let $H':=H[A_{i,1}\cup \{s_i\}]\subseteq F^*$ and $H'':=H - ((A_{i,1}\setminus \{v_{i,1}^m\})\cup \{s_i\})$.
    By Lemma~\ref{Ffaivalue}, we have $\lambda'(H')\ge Cn^{\varepsilon_0}$.
    Since $H''$ does not contain any vertex in $S$, the earlier result implies $\lambda'(H'') = \lambda(H'') \ge Cn$.
    As $H'$ and $H''$ intersects only at $\{v_{i,1}^m\}$, we have $\lambda'(H) = \lambda'(H')\lambda'(H'')n^{-1} \ge Cn^{\varepsilon_0}$.

    We now consider the case $|V(H)\cap S|\ge 2$.
    Let $L_1=\{i\in [m]: s_i\in V(H)\}$.
    For each $i\in [m]$, let $H^i := H[\{s_i\}\cup \bigcup_{j\in [t]}A_{i,j}]$ and $H' = H[A_0]$.
    We now claim that $\lambda'(H^i) \ge Cn$ for each $i \in [m]$ unless $H^i$ is the empty graph.
    If $i \notin L_1$, then $s_i \notin V(H^i)$, and $H^i$ is a subgraph of $B \setminus S$, so $\lambda'(H^i) = \lambda(H^i) \ge Cn$ by the earlier case.
    For each $i\in L_1$, since $H$ is connected, the subgraph $H^i$ must intersect each of the sets $A_{i,j}$ for all $j \in [t]$.
    Let $H^i_1:=H^i[A_{i,1}\cup \{s_i\}]$ and for $2\le j\le t  $ let $H^i_{j}:=H^i[A_{i,j}\cup \{v_{i,j-1}^m\}]$.
    Each $H^i_j$ is a subgraph of $F^*$ by construction.
    Furthermore, by the connectivity of $H^i$, we know that each $H_j^i$ contains $v_{i,j-1}^m$ if $2 \le j \le t$ and $s_i$ if $j=1$.
    Thus, by Lemma~\ref{Ffaivalue}, we have $n^{v_{H_j^i}-1}p^{e_{H_j^i}} \ge Cn^{\varepsilon_0}$ for each $ j\in [t] $.
    Therefore,
    \[
        \lambda'(H^i) = n^{v_{H^i}-1}p^{e_{H^i}} = \prod_{j\in [t]} n^{v_{H_j^i}-1}p^{e_{H_j^i}} \ge (Cn^{\varepsilon_0})^{t} \ge Cn
    \]
    where the last inequality follows by the choice of $t$.
    As $H'$ is isomorphic to a subgraph of $F$, we have $\lambda'(H') = \lambda(H') \ge Cn$.
    Let $s$ be the number of non-empty $H^i$ among $i \in [m]$.
    Observe that $H$ can be obtained by sequentially gluing each non-empty $H^i$ to $H'$ at the vertex $v_{i,t}^m$.
    Therefore, by applying Proposition~\ref{chang2phi} repeatedly, we obtain
    \[\lambda'(H) \ge n^{-s+1}\prod_{i\in [m], H^i \text{ is non-empty}} \lambda'(H^i) \ge Cn^{\varepsilon_0}.\]
    This completes the proof of Lemma~\ref{bigFfaivalue}.
\end{proof}

The following result ensures that we can find many copies of the $D$-structure in $G$.
\begin{lemma}[{\rm \cite[Corollary~2]{erdHos1983supersaturated}}]\label{supersat}
    Let $D$ be the labelled forest with $tm^2$ vertices as defined in Definition~\ref{manyF'}.
    Let $0 < \beta_{\ref{supersat}} \ll \varepsilon,\, 1/(tm^2)$.
    Then every $n$-vertex graph $G$ with at least $\varepsilon n^2$ edges contains at least $\beta_{\ref{supersat}} n^{tm^2}$ copies of $D$.
\end{lemma}

We now apply Lemmas~\ref{bigFfaivalue} and~\ref{supersat} to prove Lemma~\ref{absorblem}.

\begin{proof}[Proof of Lemma \ref{absorblem}.]

Let $\varepsilon_0=\frac{1}{m_1(F) v_F^2}$ and $C$ be the constant given in Lemma \ref{bigFfaivalue} and let $c'$ play the role of $c$ in Lemma \ref{bigFfaivalue}.
	Let $ \beta_{\ref{supersat}}$ be the constant given in Lemma \ref{supersat}.
Let $0<1/c\ll 1/c'$.
Fix constants   $1/n\ll1/C , \xi\ll \beta \ll \gamma_1 \ll \beta_{\ref{supersat}} \ll \tau  \ll \varepsilon, \varepsilon_0, 1/m$. Let $t=\lceil\frac{1}{\varepsilon_0}\rceil\le m^3$,  $t'\in \mathbb{N}$ and $t' \le \gamma_1 n$. Then $C_1=tm^2$.

Suppose $\ell =\lceil \frac{120\tau}{\gamma_1 } \rceil$ and
$p\ge c n^{-\frac{1}{m_1(F)}} $. We will expose $R:=G(n,p)$ in $\ell+1$ rounds: $R=\bigcup_{i\in [\ell+1]}R_i$, where $R_1, \ldots, R_{\ell+1}$ are independent copies of $G(n,p')$ and $(1-p')^{\ell+1}=1-p$. Note that $(1-p')^{\ell+1}\ge 1-(\ell+1)p'$, we have $p'\ge p/(\ell+1)\ge c' n^{-\frac{1}{m_1(F)}}$. Let $G^*=G\cup G(n,p)$ and $H_i=G\cup R_i$ for $i\in [\ell +1]$.

We first randomly choose a subset $V_1 \subseteq V(G)$ by selecting each vertex independently with probability $\tau/2$. By Chernoff's bound, with high probability, $\tau n/4 \leq |V_1| \leq \tau n$.
    Let $L$ be the set of vertices in $V(G)$ with degree at least $\varepsilon n/2$.
    Then $|L| \geq \varepsilon n/2$, so again by Chernoff's bound,  w.h.p. we have $|V_1 \cap L| \geq \tau \varepsilon n/4$. Furthermore, for each $v \in L$, the probability that $v$ has fewer than $\tau \varepsilon n/4$ neighbors in $V_1$ is at most $o(1/n)$, so by a union bound, w.h.p. every vertex in $L$ has at least $\tau \varepsilon n/4$ neighbors in $V_1$.
    Therefore, w.h.p. $G[V_1]$ contains at least $\tau^2 \varepsilon^2 n^2/32$ edges.
    Thus, there exists a subset $V_1\subseteq V(G)$ with $\tau n/4 \leq |V_1| \leq \tau n$ and $e_{G[V_1]}\ge \frac{\varepsilon \tau^2 n^2}{32}$.
    We now claim that this subset $V_1$ is a desired subset. The property~\ref{cond:absorbing1} comes from the choice of $V_1$.

	To prove properties \ref{cond:absorbing2} and \ref{cond:absorbing3}, we first introduce the following claim.

	\begin{claim}\label{subgembed}
		Let $X\subseteq V(G^*)$ with $  |X|\ge \tau n/4$ and  $e_{G[X]}\ge \frac{\varepsilon \tau^2 n^2}{32}$. Let  $\{S_1,\dots, S_{t'}\}$ be a  family of vertex-disjoint $m$-subsets in $V(G^*)$.
	Then for any $j\in [\ell+1]$, w.h.p.
		$H_j$ contains vertex-disjoint subsets $T_1,\dots, T_{t'}$ in $X\setminus (\bigcup_{i\in [t']}S_i)$ such that $T_i$ is an $S_i$-absorber.
	\end{claim}
	\begin{proof}
		Our goal is to find $t'$   vertex-disjoint $tm^2$-sets $T_1, \ldots, T_{t'}$ in $X\setminus (\bigcup_{i\in [t']}S_i)$ such that $T_i$ is an $S_i$-absorber in  $H_j$ for any $j\in [\ell+1]$.
		The first step is to find the candidates for $S_i$-absorbers in $X$.
		By Lemma \ref{supersat},
		there are at least $ \beta_{\ref{supersat}}  |X|^{tm^2} $  ordered $tm^2$-sets $D$ in $G[X]$. Denote by $\mathcal{F}$ the collection of  such  ordered $tm^2$-sets. Then $|\mathcal{F}|\ge  \beta_{\ref{supersat}} (\tau/4)^{tm^2} n^{tm^2} $.
		Let $\mathcal{F}_i=\mathcal{F}$ for  $ i\in  [t']$.
		The second step is to find an $S_i$-absorber  in $\mathcal{F}_i$ for any $i\in [t']$  by adding the edges of $R_j$ to fill the missing edges.
	By Lemma \ref{bigFfaivalue}, we have   $\Phi_{B\setminus S}\ge C n$ and $\Phi_{B,S}\ge Cn^{\varepsilon_0}$.
		Applying  Proposition \ref{propembed} to $R_j$, where
		$(t',\lceil\frac{1}{\varepsilon_0}\rceil,    \beta_{\ref{supersat}} (\tau/4)^{tm^2} ,  B, S, S_i)$ plays the role  of $( t,  \ell, \beta, F, W, U_i) $, w.h.p. there exist pairwise disjoint embeddings $\phi_i$ of $B$ into $R_j$ such that $S$ is mapped to some $S_i$ and $V(B)\setminus S$ is mapped to some $tm^2$-set in $\mathcal{F}_i$.
        Thus, by Lemma \ref{prop:BD_is_absorber}, we obtain $t_1$
          vertex-disjoint $tm^2$-sets $T_1, \ldots, T_{t'}$ in $X\setminus (\bigcup_{i\in [t']}S_i)$ such that $T_i \in \mathcal{F}_i$ is an $S_i$-absorber in  $H_j[X]$.
	\end{proof}

	 \textbf{Proof of property \ref{cond:absorbing2}:} For every $S\subseteq V(H)$ with $|S|\le \gamma_1 n$, there exists a family of vertex-disjoint $m$-sets $\{S_1, \ldots, S_{t'}\}$, where $t'\le \gamma_1 n$ such that $S\subseteq \bigcup_{i\in [t']} S_i $.
     Applying  Claim \ref{subgembed},  where $V_1$ plays the role of $X$,
     then w.h.p. there exists $T\subseteq V_1\setminus S$ with $|T|=(m-1)|S|$ such that $H[S\cup T]$ contains an $F$-factor.

	\textbf{Proof of property \ref{cond:absorbing3}:} For  $\{S_{i,1},\ldots,S_{i,t_i}\}$ with $1\le i\le k$ and $t_i\le \gamma_1 n$,
	we find $\{T_{i,1},\ldots,T_{i,t_i}\}$ step by step.
	 Let $t_0\in \mathbb{N}$ and $ T_{0,j}=\emptyset$ for any $j\in [t_0]$. Assume that for any $0\le x\le i-1$, $\{T_{x,1}\cup \cdots \cup T_{x,t_x}\}$ has been found.   We now find  $\{T_{i,1},\ldots,T_{i,t_i}\}$.
    Let $W=\bigcup_{i\in [k], j\in [t_i]} S_{i,j}$ and $X= V(G)\setminus\big(W\cup \bigcup_{0\le x\le i-1} (T_{x,1}\cup \cdots \cup T_{x,t_x}) \big)$.
	Then $|W\cup \bigcup_{0\le j\le i-1}(T_{j,1}\cup \cdots \cup T_{j,t_j})|\le k \gamma_1   mn +k \gamma_1   tm^2n\le 240 \tau tm^2n$  because $k\le \frac{120\tau}{\gamma_1}$.
	Since  $\tau \ll 1/m$, we get
	 $|X|   \ge n- 240 \tau tm^2n \ge \tau n$. Additionally, $e_{G[X]}\ge \varepsilon n^2- 240\tau tm^2 n^2\ge \frac{\varepsilon \tau^2  n^2}{32} $ since $t\le m^3$ and $\tau \ll \varepsilon$.
	By Claim~\ref{subgembed}, w.h.p.   $  H_{i+1}[X]$ contains  a family $\{T_{i,1}, \ldots, T_{i,t_i}\}$ such that for any $j\in [t_i]$,  $ {T}_{i,j}$ is an $ {S}_{i,j}$-absorber.
	Repeat this process for all $i\in [k]$, w.h.p. we find disjoint families  $\{T_{1,1},\ldots,T_{1,t_1}\}, \ldots$, $ \{T_{k,1},\ldots,T_{k,t_k}\}$ in $V(H)\setminus W$ such that for any $i\in [k]$ and $j\in [t_i]$, $ T_{i,j}$ is an $S_{i,j}$-absorber.

    This completes the proof.
\end{proof}

\section{Spanning bounded degree graphs in randomly perturbed graphs}

We now prove Theorem~\ref{mainbounded} using Theorem~\ref{mainfactor}.
\begin{proof}
[Proof of Theorem \ref{mainbounded}]
Let $1/n\ll 1/c \ll \gamma\ll\varepsilon,1/\Delta$.
For the given graph $F$, delete all (vertex-disjoint) copies of $K_{\Delta+1}$ and denote the resulting graph by $F'$.
Let $\ell$ be the number of copies of $K_{\Delta+1}$ in $F$.
i.e., $F=\ell K_{\Delta+1}\cup F'$.
Assume $p\geq cn^{-\frac{2}{\Delta+1}}$.

We first consider the case when $\ell(\Delta+1)\leq(1-\gamma)n$.
We expose $G(n, p)$ in two rounds: $G(n,p) = G_1 \cup G_2$, where $G_1$ and $G_2$ are independent copies of $G(n, p')$ with $(1-p')^2 = 1-p$.
Note that since $(1-p')^2>1-2p'$, we have $p'>p/2$.
We first find $\ell K_{\Delta+1}$ in $G_1$, which is guaranteed w.h.p.\ by Lemmas~\ref{otherFcover} and~\ref{otherfaiF}.
Let $R$ be the set of vertices of $\ell K_{\Delta+1}$ in $G_1$.
Note that $G_1$ and $G_2$ are independent so $R$ is independent of $G_2$.
By apply Lemma~\ref{FKNP7.6} to $F'$, the graph $G_2[V(G) \setminus R]$ contains a copy of $F'$ w.h.p.
This yields a copy of $F = \ell K_{\Delta+1} \cup F'$ in $G(n,p)$.

We now consider the case when $\ell(\Delta+1)>(1-\gamma)n$.
Again, expose $G(n, p)$ as $G(n,p) = G_1 \cup G_2$, where $G_1, G_2 \sim G(n, p')$ and $p' > p/2$ as before.
The first step is to find $F'$ in $G_1$.
By applying Lemma~\ref{FKNP7.6} to $F'$, we obtain that w.h.p.~$G_1$ contains a copy of $F'$, denoted by $Q$.
As $V(Q)$ is fixed before revealing $G_2$, it is independent of $G_2$.
Note that the $G\setminus V(Q)$ still has at least $\varepsilon n^2/2$ edges.
Now apply Theorem~\ref{mainfactor} to $G \cup G_2[V \setminus V(Q)]$ to find a $K_{\Delta+1}$-factor w.h.p..
Therefore, combining the copy of $F'$ in $Q$ and the $K_{\Delta+1}$-factor in $V(G) \setminus V(Q)$, we conclude that  w.h.p. $G \cup G(n,p)$ contains a copy of $F$.
\end{proof}

\section{Hamilton $\ell$-cycles in randomly pertured $k$-graphs}

Given a hypergraph $H=(V,E)$, we define a \emph{walk} in $H$ to be an alternating sequence $v_1,e_1,v_2,\ldots,$
$e_s,v_{s+1}$ of vertices and edges of $H$ such that the following folds:
\begin{itemize}
  \item $v_j\in V$ for each $j\in[s+1]$,
  \item $e_j\in E$ for each $j\in[s]$,
  \item $\{v_j,v_{j+1}\}\subseteq e_j$ for each $j\in[s]$.
\end{itemize}
Such a walk from $v_1$ to $v_{s+1}$ will be referred as a $(v_1,v_{s+1})$-walk.
A \emph{path} is a walk with the additional restrictions that the $s+1$ vertices are all distinct and the $s$ edges are all
distinct.
Two vertices $v,w\in V$ are said to be \emph{connected} in $H$ if there exists a $(v,w)$-path in $H$; otherwise, $v$ and $w$ are separated from each other.
A hypergraph $H$ is \emph{connected} if every pair of vertices $v,w \in V$ is connected in $H$; otherwise $H$ is \emph{disconnected}.
A \emph{component} of a hypergraph $H$ is a maximal connected induced subhypergraph of $H$.
The following observation allows us to focus on connected subgraphs.
\begin{lemma}\label{lem:components}
    Let $H$ be a $k$-graph with $a$ components $H_1, \ldots, H_a$.
    Then $d(H) \leq \max\{d(H_1), \ldots, d(H_a)\}$.
\end{lemma}
\begin{proof}
    By definition, we  see that
\begin{equation}
\nonumber
\begin{split}
d(H)&=\frac{e_{H}}{v_{H}-1}=\frac{\sum_{i\in[a]}e_{H_i}}{\sum_{i\in[a]}v_{H_i}-1}\leq\frac{\sum_{i\in[a]}e_{H_i}}{\sum_{i\in[a]}v_{H_i}-a}\\
&\leq\max\left\{\frac{e_{H_1}}{v_{H_1}-1},\ldots,\frac{e_{H_a}}{v_{H_a}-1}\right\}\\
& = \max\{d(H_1), \ldots, d(H_a)\}.\qedhere
\end{split}
\end{equation}

\end{proof}

Let $\ell \ge 2$ and $\text{HC}_{\ell}$ be the labelled Hamilton $\ell$-cycle on the vertex set $x_1, \ldots, x_n$, with the obvious cyclic order and on the edge set $e_1,\ldots,e_{\frac{n}{k-\ell}}$.
Fix an integer $M\in(k-\ell)\mathbb{N}$ with $M \geq 2k^2$, which will be specified later.
For each $j \in [\lfloor\frac{n}{M}\rfloor]$, we define $P_j$ to be the induced subgraph of $\text{HC}_{\ell}$ on the vertex set $\{x_{jM+1},\ldots,x_{jM+k+(\lceil\frac{k}{k-\ell}\rceil-1)(k-\ell)}\}$.
Note that $P_j$ is an $\ell$-path with $\lceil\frac{k}{k-\ell}\rceil$ edges.
Let $P = \bigcup_{j \in [\lfloor\frac{n}{M}\rfloor]} P_j$ and let $\overline{\text{HC}_{\ell}}$ be the subgraph of $\text{HC}_{\ell}$ obtained by removing $E(P)$.
Then we note that each component of $\overline{\text{HC}_{\ell}}$ is an $\ell$-path with at most $2M$ vertices.

\begin{lemma}
\label{m1HC}
If $1/n\ll1/M\ll1/k,1/\ell$, then $m_1(\overline{\text{HC}_{\ell}})\leq\frac{1}{k-\ell+\frac{1}{2M}}$.
\end{lemma}
\begin{proof}
Let $H'$ be an induced subgraph of $\overline{\text{HC}_{\ell}}$ on $v$ vertices. By Lemma~\ref{lem:components}, we may assume that $H'$ is connected.
As each component of $\overline{\text{HC}_{\ell}}$ is a subgraph an $\ell$-path with at most $2M$ vertices, $H'$ is a subgraph of an $\ell$-path with at most $2M$ vertices.
Thus, we can observe that
\begin{equation}
\nonumber
\begin{split}
d(H') \leq\frac{2M-\ell}{(k-\ell)(2M-1)}\leq\frac{1}{k-\ell+\frac{1}{2M}},
\end{split}
\end{equation}
where the last inequality holds since $\frac{2M-\ell}{(k-\ell)(2M-1)}=\frac{1}{(k-\ell)(1+\frac{2M-1}{2M-\ell})} $ and $\frac{(k-\ell)(2M-1)}{2M-\ell}\geq\frac{1}{2M}$.
As $H'$ is arbitrary, we conclude that $m_1(\overline{\text{HC}_{\ell}})\leq\frac{1}{k-\ell+\frac{1}{2M}}$.
\end{proof}

Moreover, we also need the following supersaturation lemma, which asserts that every $k$-graph with positive density contains many $\ell$-paths of constant length.
\begin{lemma}
[{\rm \cite[Corollary~2]{erdHos1983supersaturated}}]
\label{su}
  Let $k' = k+(\lceil\frac{k}{k-\ell}\rceil-1)(k-\ell)$.
  For any $\varepsilon>0$, there exists $\varepsilon'>0$ such that if $H$ is an $n$-vertex $k$-graph with at least $\varepsilon n^k$ edges, then it contains at least $\varepsilon'n^{k'}$ labeled copies of $\ell$-paths with $k'$ vertices.
\end{lemma}
\begin{proof}[Proof of Theorem~\ref{hypergraph}]
For given $\varepsilon, r$, we choose $M\in(k-\ell)\mathbb{N}$ and $\eta>0$ such that \[1/n\ll\eta\ll1/M\ll\varepsilon'\ll\varepsilon, 1/\ell,1/k.\]
Let $H$ be an $n$-vertex $k$-graph with $e(H) \ge \varepsilon n^k$.
Let $\overline{\text{HC}_{\ell}}$ be the subgraph of a Hamilton $\ell$-cycle $\text{HC}_{\ell}$ defined above with the choice of $M$.
Let $k' = k+(\lceil\frac{k}{k-\ell}\rceil-1)(k-\ell)$.

We now define a $2$-edge-colored complete $k$-graph $K$ on $V(H)$, where an edge is colored blue if and only if it belongs to $E(H)$.
Let $\Phi$ be the family of embeddings of $\text{HC}_{\ell}$ into $K$ such that each edge in $P$ is mapped to a blue edge.
Note that in the randomly perturbed model $H \cup G^{(k)}(n,p)$, if there exists an embedding $\phi \in \Phi$ such that all edges of $\phi(\overline{\text{HC}_{\ell}})$ are present in $G^{(k)}(n,p)$, then $H \cup G^{(k)}(n,p)$ contains a copy of $\text{HC}_{\ell}$.
We will show that the family $\Phi$ admits a $(2^{k'+1}/\varepsilon')n^{-1}$-vertex spread distribution.

Let $Q$ be an $\ell$-path with $k'$ vertices and let $V(Q) = \{w_1, \ldots, w_{k'}\}$ where the vertices are labeled by a natural ordering. By Lemma~\ref{su}, we obtain that there are at least $\varepsilon'n^{k'}$ labeled copies of $Q$ in $H$.
Let $L_1$ be the set of vertices $y_1$ in $V(K)$ such that there are at least $\varepsilon' n^{k'-1}/2$ embeddings of $Q$ in $H$ that maps $w_1$ to $y_1$.
Then $|L_1| \geq \varepsilon' n/2$ as otherwise the number of possible embeddings of $Q$ in $H$ would be less than $|L_1| \times n^{k'-1} + n \times \varepsilon' n^{k-1}/2 \leq \varepsilon' n^{k'}$, which is a contradiction.
For each $y_1 \in L_1$, let $L_2(y_1)$ be the set of vertices $y_2$ in $V(K)$ such that there are at least $\varepsilon' n^{k'-2}/4$ embeddings of $Q$ in $H$ that maps $w_1$ to $y_1$ and $w_2$ to $y_2$.
Then by the same argument, we have $|L_2(y_1)| \geq \varepsilon' n/4$ for each $y_1 \in L_1$.
By applying the same construction iteratively, we obtain sets $L_i(y_1, \ldots, y_{i-1}) \subseteq V(K)$ for each $i \in [k']$ such that the following holds.
\begin{itemize}
    \item For each $i \in [k']$ and $y_1 \in L_1, y_2 \in L_2(y_1), \ldots, y_{i-1} \in L_{i-1}(y_1, \ldots, y_{i-2})$, we have $|L_i(y_1, \ldots, y_{i-1})| \geq \varepsilon' n/2^i$.
    \item For each $y_1 \in L_1, y_2 \in L_2(y_1), \ldots, y_{k'} \in L_{k'}(y_1, \ldots, y_{k'-1})$, the map $w_i \mapsto y_i$ for $i \in [k']$ is an embedding of $Q$ into $H$.
\end{itemize}

We now describe a randomized embedding process from $\{x_1,\ldots, x_n\}$ to $V(K)$.
We first embed $x_{jM+1},\ldots,x_{jM+k'}$ for each $j\in[\lfloor\frac{n}{M}\rfloor]$.
For each $j \in [\lfloor\frac{n}{M}\rfloor]$, assume that $x_{j'M+1},\ldots,x_{j'M+k'}$ are already embedded for $j' < j$ and let $X$ be the embedded image of these vertices.
We first embed $x_{jM+1}$ uniformly at random from $L_1 \setminus X$.
Since at most $k'\lfloor\frac{n}{M}\rfloor$ vertices are already assigned, there are at least $\varepsilon'n/2 - k'\lfloor\frac{n}{M}\rfloor\geq \varepsilon' n/4 \geq \varepsilon n/2^{k'+1}$ choices for the image of $x_{jM+1}$.
Suppose that $x_{jM+1}$ is embedded into $y_1 \in L_1 \setminus X$.
We embed $x_{jM+2}$ uniformly at random from $L_2(y_1) \setminus X$.
By the similar argument, we have at least $|L_2(y_1) \setminus X| \geq \varepsilon' n/4 - k'\frac{n}{M} \geq \varepsilon n/2^{k'+1}$ choices for the image of $x_{jM+2}$.
We apply this process iteratively: For each $i \in [k']$, conditioning on $x_{jM+1}, x_{jM+2}, \ldots, x_{jM+i-1}$ being embedded into $y_1, y_2, \ldots, y_{i-1}$ respectively, we embed $x_{jM+i}$ uniformly at random from $L_i(y_1, \ldots, y_{i-1}) \setminus X$.
Then we have at least $\varepsilon n/2^{k'+1}$ choices for the image of $x_{jM+i}$.

After embedding all the vertices in $V(P)$, the vertices not in $V(P)$ are mapped uniformly at random from the remaining vertices.
This defines a probability distribution on the family $\Phi$.
We now claim that this is a $(2^{k'+1}/\varepsilon')n^{-1}$-vertex-spread distribution on $\Phi$.

Let $C=2^{k'+1}/\varepsilon'$.
Recall that in our random process of generating the embeddings of $\text{HC}$ into $K$, we first choose the images of the vertices in $V(P)$ and then assign images to the remaining vertices.
For any $s\in[n]$ and any sequences $y_1,\ldots,y_s\in V(\overline{\text{HC}})$, $z_1,\ldots,z_s\in V(K)$, let $t:=|\{y_1,\ldots,y_s\}\cap V(P)|$.
If $s=t$, then $\mu(\{\psi:\psi(y_i)=z_i, \forall \ i\in[s]\})\leq\frac{1}{(\varepsilon' n/2^{k'+1})^s}\leq\left(\frac{C}{n}\right)^s$.
If $s\geq t+1$, then we have
\begin{equation}\nonumber
\begin{split}
\mu\left(\{\psi:\psi(y_i)=z_i, \forall \ i\in[s]\}\right)&\leq\frac{1}{(\varepsilon' n/2^{k'+1})^t\cdot(n-k'\lfloor\frac{n}{M}\rfloor)\cdots(n-k'\lfloor\frac{n}{M}\rfloor-(s-t-1))}\\
&=\left(\frac{C}{n}\right)^t\frac{(n-k'\lfloor\frac{n}{M}\rfloor-(s-t))!}{(n-k'\lfloor\frac{n}{M}\rfloor)!}\\
&\leq\left(\frac{C}{n}\right)^t\left(\frac{e}{n-k'\lfloor\frac{n}{M}\rfloor}\right)^{s-t} \le \left(\frac{C}{n}\right)^{s}
\end{split}
\end{equation}
where third inequality uses $1/M\ll 1/\ell,1/k$.
Therefore, $\Phi$ admits a $(C/n)$-vertex spread distribution.

We now consider each $\phi \in \Phi$ as an embedding of $\overline{\text{HC}_{\ell}}$ into $K$.
By Proposition~\ref{kelly}, there exists a constant $C' =C'(\varepsilon')$ such that there is a $(C'/n^{1/m_1(\overline{\text{HC}_{\ell}})})$-spread distribution on subgraphs of $K$ isomorphic to $\overline{\text{HC}_{\ell}}$, which admits a labeling such that all edges of $P$ are blue edges (edges of $H$).

By Lemma~\ref{m1HC}, we obtain that $m_1(\overline{\text{HC}_{\ell}})\leq \frac{1}{k-\ell+1/(2M)}$.
Let $p = n^{-(k-\ell)-\eta}$. Then since $1/n\ll\eta\ll1/M \ll1/k,1/\ell$, we have $p \geq n^{-(k-\ell)-1/(2M)+\eta}\geq n^{-\frac{1}{m_1(\overline{\text{HC}_{\ell}})}}\log n$.
Therefore, by Proposition~\ref{FKNP}, with high probability the random graph $G^{(k)}(n,p)$ contains a copy of $\overline{\text{HC}_{\ell}}$ such that all edges in $E$ appear in $H$.
This yields a copy of Hamilton $\ell$-cycle in $H\cup G^{(k)}(n,p)$.
\end{proof}

\section{Powers of tight Hamilton cycles in randomly perturbed $k$-graphs}

We now turn to the proof of Theorem~\ref{mainpower}.
The proof is similar to that of Theorem~\ref{hypergraph} but more involved.
Indeed, we may not be able to disconnect the power of tight cycle by removing small sets of edges of \textit{zero Tur\'an density}.
So we have to compute the $1$-density of the target graph more carefully.

Let $r \ge 2$ and $\text{HC}^r$ be the labelled $r$-th power of a tight Hamilton cycle on the vertex set $x_1, \ldots, x_n$, with the obvious cyclic order.
Fix a sufficiently large integer $M$, which will be specified later.
We define a set of edges $E:=\{\{x_{ jM+1},x_{jM+1},\ldots,x_{jM+k}\}: j\in[\lfloor n/M\rfloor]\}$.
Let $\overline{\text{HC}^r}$ be the subgraph of $\text{HC}^r$ obtained by removing the edges of $E$.
We first count the number of edges of induced subgraphs of $\overline{\text{HC}^r}$.
\begin{proposition}
  \label{factnumber}
Let $P$ be a connected induced subgraph of $\overline{\text{HC}^r}$ on $w$ vertices with $w \geq k+r-1$.
Then $e_{P}\le w\binom{k+r-2}{k-1}-w/M+2$.
Furthermore, if $P$ is a subgraph of an $r$-th power of a tight path, then $e_{P}\leq\binom{k+r-1}{k}+(w-(k+r-1))\binom{k+r-2}{k-1}-w/M+2\le w\binom{k+r-2}{k-1}-w/M$.
\end{proposition}
\begin{proof}
We first consider the case when $P$ is a subgraph of $\overline{\text{HC}^r}$ and also a subgraph of an $r$-th power of a tight path.
Let $P'$ be the induced subgraph of $\text{HC}^r$ on $V(P)$.
Suppose that $P$ begins at $x_{q+1}$ and ends at $x_{q + w'}$ where the indices are considered in modulo $n$.
Note that $w' \ge w$.

We now aim to estimate $e(P)$.
Assume that $V(P)$ consists of $\ell$ intervals $V_1, \ldots, V_{\ell}$ along the cyclic order, and let $d_p$ be the number of vertices between intervals $V_p$ and $V_{p+1}$ for each $p \in [\ell - 1]$.
Note that $\sum_{p \in [\ell - 1]} d_p = w' - w$, and $d_p \le k+r-1$ since $P$ is connected.

We now count $e_{P'}$ by sequentially adding the vertices of $P'$ in their cyclic order.
Let $v_1, \ldots, v_{w}$ be the vertices of $P$ in cyclic order and we say an edge $e$ a \emph{backward edge} of $v_j$ if $e$ is incident to $v_j$ and the other vertices incident to $e$ are those that appear before $v_j$ in the cyclic order.
Let $h_j=\binom{j-1}{k-1}$ if $j\le k+r-1$ and $h_j=\binom{k+r-2}{k-1}$ otherwise.
Observe that at $v_j$ has most $h_j$ backward edges.
Moreover, if $v_j$ is the first vertex of an interval $V_{p+1}$ for some $p \in [\ell - 1]$, then it has at most $h_j - d_p$ backward edges due to the gap of $d_p$ vertices between the interval.
Thus, the total number of edges in $P'$ satisfies $e_{P'}\leq \sum_{j\in [w]} h_j - \sum_{p\in [\ell-1]}d_p=\binom{k+r-1}{k}+(w-(k+r-1))\binom{k+r-2}{k-1}-(w'-w)$.

Observe that in any subgraph of $\text{HC}^r$ on $w$ consecutive vertices, there are at least $\lceil w / M \rceil - 2$ edges of $E$.
If an edge $e \in E$ lies in the interval $\{x_{q+1},\ldots, x_{q + w'}\}$ but is not in $P'$, then at least one of its vertices must lie outside $V(P')$.
As there are $w' - w$ vertices in $\{x_{q+1},\ldots, x_{q+w'}\} \setminus V(P')$, we obtain $|E(P'_i)\cap E|\geq \lceil w/M\rceil-2-(w'-w)$.
We recall that $E(P) = E(P') \setminus E$.
Therefore, we get
\[e_{P}\leq \binom{k+r-1}{k}+(w-(k+r-1))\binom{k+r-2}{k-1}-\left(w/M-2\right)\leq w\binom{k+r-2}{k-1}-w/M\]
where the last inequality is from $2\leq\frac{k-1}{k}(k+r-1)\binom{k+r-2}{k-1}$.
This proves the second bound in the proposition.

We now consider the case when $P$ is a subgraph of $\overline{\text{HC}^r}$ but is not a subgraph of an $r$-th power of a tight path.
In this case, $P$ is connected and if we order the vertices of $V(P)$ with respect to the cyclic order of $\overline{\text{HC}^r}$, then the cyclic distance between any two consecutive vertices is at most $k+r-1$.

We count $e_P$ as before by adding vertices sequentially following the cyclic order.
As $P$ may not have clear ends, we slightly modify the definition of backward edges.
For $v = x_i \in V(P)$ and an edge $e$ incident to $v$, we say $e$ is an \emph{backward edge} of $v$ if $e$ is contained in $\{x_{i-(k+r-1)},\ldots,x_{i-1}, x_i\}$.
Then the difference of the proof is that even the first $k+r-1$ vertices may have up to $\binom{k+r-2}{k-1}$ backward edges.
Let $V_1, \ldots, V_\ell$ be the intervals in $V(P)$ along the cycle, and let $V_{\ell + 1} = V_1$.
Let $d_p$ be the number of vertices between intervals $V_p$ and $V_{p+1}$ for each $p \in [\ell]$.
Then the first vertex of interval $V_{p+1},p\in[\ell]$ has at most $\binom{k+r-2}{k-1}-d_p$ backward neighbors.
Thus, we obtain that $e_{P'}\leq w\binom{k+r-2}{k-1} -\sum_{p\in[\ell]}d_p = w\binom{k+r-2}{k-1}-(w'-w)$ where $P'$ is an induced subgraph of $\text{HC}^r$ on $V(P)$.
By the same argument as above, we have $e_P\leq w\binom{k+r-2}{k-1}-w/M+2$.
\end{proof}

This proposition allows us to bound $m_1(\overline{{\emph HC}^r})$.
\begin{lemma}
\label{m1HCr}
If $1/n \ll 1/M\ll1/r,1/k$, then $m_1(\overline{{\emph HC}^r})\leq \binom{k+r-2}{k-1}-1/(2M)$.
\end{lemma}
\begin{proof}
Let $H'$ be a connected induced subgraph of $\overline{\text{HC}^r}$ on $v$ vertices.
If $v\leq k+r-1$, then
\begin{equation}\nonumber
\begin{split}
d(H')&\le d(K_v^{(k)}) =\frac{\binom{v}{k}}{v-1}=\frac{v(v-2)\cdots(v-k+1)}{k!}\\
&\leq \frac{(k+r-1)(k+r-3)\cdots r}{k!}\leq \frac{(k+r-2)(k+r-3)\cdots r}{(k-1)!}-\frac{1}{2M}\\
&=\binom{k+r-2}{k-1}-1/(2M),
\end{split}
\end{equation}
where the third inequality holds since $k+r-1\leq k(k+r-2)-\frac{k!}{2M(k+r-3)\cdots r}$ for $1/M \ll 1/k, 1/r$.

If $k+r\leq v\leq M\binom{k+r-2}{k-1}+1$, then $H'$ is a subgraph of an \emph{r}-th power of a tight path. Thus, we have
\begin{equation}\nonumber
\begin{split}
d(H')&=\frac{e_{H'}}{v-1}\leq\frac{\binom{k+r-1}{k}+(v-(k+r-1))\binom{k+r-2}{k-1}}{v-1}\\
&=\binom{k+r-2}{k-1}-\frac{(k+r-2)\binom{k+r-2}{k-1}-\binom{k+r-1}{k}}{v-1}\\
&\leq\binom{k+r-2}{k-1}-\frac{(k+r-2)\binom{k+r-2}{k-1}-\binom{k+r-1}{k}}{M\binom{k+r-2}{k-1}}\\
&\leq\binom{k+r-2}{k-1}-\frac{k+r-2-\frac{k+r-1}{k}}{M}\\
&\leq\binom{k+r-2}{k-1}-\frac{1}{2M},
\end{split}
\end{equation}
where the last inequality holds since $k+\frac{1}{k}+r-3-\frac{r}{k}\geq2+\frac{1}{2}-3+\frac{r}{2}\geq\frac{1}{2}$.

If $\binom{k+r-2}{k-1}+1\leq v<n/(k+r-1)$,
then there exist two vertices of $V(H')$ whose distance with respect to the cyclic order is larger than $k+r-1$, and no vertex of $H'$ lies between them in the cyclic order.
Thus, $H'$ is a subgraph of the $r$-th power of a tight path.

By Proposition~\ref{factnumber}, we obtain that $e_{H'}\leq\binom{k+r-1}{k}+(v-(k+r-1))\binom{k+r-2}{k-1}-v/M+2$.
Therefore, we have
\begin{equation}\nonumber
\begin{split}
d(H')&=\frac{e_{H'}}{v-1}\leq\frac{\binom{k+r-1}{k}+(v-(k+r-1))\binom{k+r-2}{k-1}-v/M+2}{v-1}\\
&=\binom{k+r-2}{k-1}-\frac{(k+r-2-\frac{k+r-1}{k})\binom{k+r-2}{k-1}+v/M-2}{v-1}\\
&\leq\binom{k+r-2}{k-1}-\frac{1}{2M},
\end{split}
\end{equation}
where the last inequality holds since $v\geq M\binom{k+r-2}{k-1}+1$, $k+r-2-\frac{k+r-1}{k}\geq\frac{1}{2}$ (which is the same as the second case) and $\binom{k+r-2}{k-1}\geq2$.

If $v\geq n/(k+r-1)$, then by
Proposition~\ref{factnumber}, we have $e_{H'}\leq v\binom{k+r-2}{k-1}-v/M+2$. Thus,
\[
d(H')=\frac{e_{H'}}{v-1}\leq\frac{v\binom{k+r-2}{k-1}-v/M+2}{v-1}\le \binom{k+r-2}{k-1}-\frac{v/M-\binom{k+r-2}{k-1}-2}{v-1}\leq \binom{k+r-2}{k-1}-\frac{1}{2M},
\]
where the last inequality uses the facts that $v\geq n/(k+r-1)$ and $n$ is sufficiently large.

In all cases, we have $d(H') \le \binom{k+r-2}{k-1}-\frac{1}{2M}$. Since $H'$ is an arbitrary connected subgraph, by Lemma~\ref{lem:components}, this completes the proof.
\end{proof}

\begin{proof}[Proof of Theorem~\ref{mainpower}]
For given $\varepsilon, r$, we choose $M$ and $\eta>0$ such that \[1/n\ll\eta\ll1/M\ll\varepsilon, 1/r.\]
Let $H$ be an $n$-vertex $k$-graph with $e(H) \ge \varepsilon n^k$.
Let $\overline{\text{HC}^r}$ be the subgraph of an $r$-th power of a tight Hamilton cycle $\text{HC}^r$ defined above with the choice of $M$.

We now define a 2-edge-colored complete $k$-graph $K$ on $V(H)$, where an edge is colored blue if and only if it belongs to $E(H)$.
Let $\Phi$ be the family of embeddings of $\text{HC}^r$ into $K$ such that each edge in $E$ is mapped to a blue edge.
Note that in the randomly perturbed model $H \cup G^{(k)}(n,p)$, if there exists an embedding $\phi \in \Phi$ such that all edges of $\phi(\overline{\text{HC}^r})$ are present in $G^{(k)}(n,p)$, then $H \cup G^{(k)}(n,p)$ contains a copy of $\text{HC}^r$.
We will show that the family $\Phi$ admits a $(2^{k+1}/\varepsilon)n^{-1}$-vertex spread distribution.

Let $L_1=\{y_1\in V(H):d_H(y_1)\geq\varepsilon n^{k-1}/2\}$.
Since $\varepsilon n^k\leq e_H\leq|L_1|\cdot n^{k-1}+\varepsilon n^{k-1}/2\cdot n$, we can obtain $|L_1|\geq \varepsilon n/2$.
For each $y_1 \in L_1$, we define $L_2(y_1)=\{y_2\in V(H): d_H(y_1, y_2)\geq\varepsilon n^{k-2}/4\}$.
Similarly, $\varepsilon n^{k-1}/2 \leq d_H(y_1) \leq |L_2(y_1)| \cdot n^{k-2} + \varepsilon n^{k-2}/4 \cdot n$. Thus, we can obtain $|L_2(y_1)|\geq \varepsilon n/4$.
We apply this process iteratively to define sets $L_3(y_1, y_2), \ldots, L_k(y_1, \ldots, y_{k-1})$ such that $|L_i(y_1, \ldots, y_{i-1})|\geq \varepsilon n/2^{i}$ for all $i\in[k]$ and $y_1 \in L_1, y_2 \in L_2(y_1), \ldots, y_r \in L_k(y_1, \ldots, y_{k-1})$ implies $\{y_1, \ldots, y_k\}$ is an edge of $H$.

We now describe a randomized embedding process from $\{x_1,\dots, x_n\}$ to $V(K)$.
For each $j \in [\lfloor n/M \rfloor]$, assume that $x_{j'M+1},\ldots,x_{j'M+k}$ are already embedded for $j'<j$ and let $X$ be the embedded image.
We first choose the image of $x_{jM+1}$ from $L_1 \setminus X$ uniformly at random. Since at most $kn/M$ vertices are already assigned, there are at least $|L_1|-kn/M \geq \varepsilon n/4 \geq \varepsilon n/2^{k+1}$ choices for the image of $x_{jM+1}$.
Conditioning on $x_{jM+1}$ being embedded into $y_1$, we embed $x_{jM+2}$ into $L_2(y_1) \setminus X$ uniformly at random.
By the same computation as before, at least $|L_2(y_1) \setminus X| \geq \varepsilon n/8 \geq \varepsilon n/2^{k+1}$ possible choices for the image of $x_{jM+2}$.
We iteratively apply this process for $x_{jM+3},\ldots,x_{jM+k}$: Conditioning on $x_{jM+1}, \ldots, x_{jM+i}$ being embedded into $y_1, \ldots, y_{i}$, we embed $x_{jM+i+1}$ into $L_{i+1}(y_1, \ldots, y_{i}) \setminus X$ uniformly at random.
Then there are at least $|L_{i+1}(y_1, \ldots, y_{i}) \setminus X| \geq \varepsilon n/2^{i+2} \geq \varepsilon n/2^{k+1}$ possible choices for the image of $x_{jM+i+1}$ for $i\in[0,k-1]$.

After we embed all the vertices in $V(E)$, we embed the vertices in $V(\overline{\text{HC}^r}) \setminus V(E)$ by uniformly at random from the remaining vertices. This defines a probability distribution on the family $\Phi$.
We now claim that this is a $(2^{k+1}/\varepsilon)n^{-1}$-vertex-spread distribution on $\Phi$.

Let $C=2^{k+1}/\varepsilon$.
Recall that in our random process of generating the embeddings of $\text{HC}^r$ into $K$, we first choose the images of the vertices in $V(E)$ and then assign images to the remaining vertices.
For any $s\in[n]$ and any sequences $y_1,\ldots,y_s\in V(\overline{\text{HC}^r})$, $z_1,\ldots,z_s\in V(K_n^{(k)})$, let $t:=|\{y_1,\ldots,y_s\}\cap\bigcup_{j\in[\lfloor n/M\rfloor]}\{x_{jM+1},x_{jM +2},\ldots,x_{jM+k}\}|$.
If $s=t$, then $\mu(\{\psi:\psi(y_i)=z_i, \forall \ i\in[s]\})\leq\frac{1}{(\varepsilon n/2^{k+1})^s}\leq\left(\frac{C}{n}\right)^s$.
If $s\geq t+1$, then we have
\begin{equation}\nonumber
\begin{split}
\mu\left(\{\psi:\psi(y_i)=z_i, \forall \ i\in[s]\}\right)&\leq\frac{1}{(\varepsilon n/2^{k+1})^t\cdot(n-k\lfloor n/M\rfloor)\cdots(n-k\lfloor n/M\rfloor-(s-t-1))}\\
&=\left(\frac{C}{n}\right)^t\frac{(n-k\lfloor n/M\rfloor-(s-t))!}{(n-k\lfloor n/M\rfloor)!}\\
&\leq\left(\frac{C}{n}\right)^t\left(\frac{e}{n-k\lfloor n/M\rfloor}\right)^{s-t} \le \left(\frac{C}{n}\right)^{s}
\end{split}
\end{equation}
where third inequality uses $1/M\ll 1/r$.
Therefore, $\Phi$ admits a $(C/n)$-vertex spread distribution.

By Proposition~\ref{kelly}, there exists a constant $C' =C'(\varepsilon)$ such that there is a $(C'/n^{1/m_1(\overline{\text{HC}^r})})$-spread distribution on subgraphs of $K$ isomorphic to $\overline{\text{HC}^r}$, which admits a labeling such that all edges of $E$ are blue edges (edges of $G$).

By Lemma~\ref{m1HCr}, we obtain that $m_1(\overline{\text{HC}^r})\leq \binom{k+r-2}{k-1}-\frac{1}{2M}$.
Let $p = n^{-\binom{k+r-2}{k-1}-\eta}$. Then since $1/n\ll\eta\ll1/M \ll1/r,1/k$, we have
\[
p \geq n^{-1/\left(\binom{k+r-2}{k-1}-\frac{1}{2M}\right)+\eta}\geq n^{-1/{m_1(\overline{\text{HC}^r})}}\log n
\]
by using $1/n\ll\eta\ll1/M \ll1/r,1/k$.
Therefore, by Proposition~\ref{FKNP}, w.h.p. the random graph $G^{(k)}(n,p)$ contains a copy of $\overline{\text{HC}^r}$ such that all edges in $E$ appear in $H$.
This yields a copy of the $r$-th power of a tight Hamilton cycle in $H\cup G^{(k)}(n,p)$.
\end{proof}

\section{Concluding Remarks}

In this paper, we studied a relaxation of the randomly perturbed graph model, namely, replacing the minimum degree condition by a density condition, when the minimum degree assumption is a vanishing one.
As mentioned above, such a weakening is not possible for $p=o(\log n/n)$, as then $G\cup G(n,p)$ may contain isolated vertices and thus does not contain any meaningful spanning subgraphs.
This observation extends to $k$-uniform hypergraphs, saying that such a weakening might be possible only if $p=\omega(\log n/n^{k-1})$, and thus excludes perfect matching and loose Hamilton cycle from the discussions.
On the other hand, the analogue of Theorem~\ref{thm:BTW} in uniform hypergraphs is currently unknown (even for \textit{clique}-factors), so a hypergraph version of Theorem~\ref{mainfactor} remains open.

Note that the results~\cite{aigner2023cycle, HPH,MR4052848,HPER,NT} on randomly perturbed graphs assume non-trivial minimum degree conditions, that is, $\delta(G)\ge \alpha n$ for some fixed $\alpha\in (0,1)$.
By the reasoning in the introduction, replacing the minimum degree conditions in them by density conditions would not yield interesting results, as a dense graph with constant density may still contain a large portion of isolated vertices.

Finally, we used the (vertex)-spread method in the proof of Theorems~\ref{hypergraph} and
~\ref{mainpower}.
In a follow-up paper, we shall develop further on this method and explore the perturbation threshold for the family of $d$-degenerate graphs.

\section*{Acknowledgements}
JH was partially supported by the Natural Science Foundation of China (12371341).
SI was supported by the National Research Foundation of Korea (NRF) grant funded by the Korea government(MSIT) No. RS-2023-00210430, and supported by the Institute for Basic Science (IBS-R029-C4).
 JZ was supported by the China Postdoctoral Science Foundation (No. 2024M764113).
\bibliographystyle{plain}
\bibliography{ref}

\begin{thebibliography}{10}

\bibitem{aigner2023cycle}
E.~Aigner-Horev, D.~Hefetz, and M.~Krivelevich.
\newblock Cycle lengths in randomly perturbed graphs.
\newblock {\em Random Structures \& Algorithms}, 63(4):867--884, 2023.

\bibitem{HPH}
S.~Antoniuk, A.~Dudek, C.~Reiher, A.~Ruci\'nski, and M.~Schacht.
\newblock High powers of {H}amiltonian cycles in randomly augmented graphs.
\newblock {\em J. Graph Theory}, 98(2):255--284, 2021.

\bibitem{MR3922775}
J.~Balogh, A.~Treglown, and A.~Z. Wagner.
\newblock Tilings in randomly perturbed dense graphs.
\newblock {\em Combin. Probab. Comput.}, 28(2):159--176, 2019.

\bibitem{MR4025389}
W.~Bedenknecht, J.~Han, Y.~Kohayakawa, and G.~O. Mota.
\newblock Powers of tight {H}amilton cycles in randomly perturbed hypergraphs.
\newblock {\em Random Structures \& Algorithms}, 55(4):795--807, 2019.

\bibitem{Bennett2017AddingRE}
P.~Bennett, A.~Dudek, and A.~Frieze.
\newblock Adding random edges to create the square of a {H}amilton cycle.
\newblock {\em arXiv:1710.02716v1}, 2017.

\bibitem{MR1943857}
T.~Bohman, A.~Frieze, and R.~Martin.
\newblock How many random edges make a dense graph {H}amiltonian?
\newblock {\em Random Structures \& Algorithms}, 22(1):33--42, 2003.

\bibitem{bondy1976graph}
J.~A. Bondy and U.~S.~R. Murty.
\newblock {\em Graph theory with applications}, volume 290.
\newblock Macmillan London, 1976.

\bibitem{MR4025392}
J.~B\"ottcher, J.~Han, Y.~Kohayakawa, R.~Montgomery, O.~Parczyk, and Y.~Person.
\newblock Universality for bounded degree spanning trees in randomly perturbed
  graphs.
\newblock {\em Random Structures \& Algorithms}, 55(4):854--864, 2019.

\bibitem{MR4130332}
J.~B\"ottcher, R.~Montgomery, O.~Parczyk, and Y.~Person.
\newblock Embedding spanning bounded degree graphs in randomly perturbed
  graphs.
\newblock {\em Mathematika}, 66(2):422--447, 2020.

\bibitem{chang2022factors}
Y.~Chang, J.~Han, Y.~Kohayakawa, P.~Morris, and G.~O. Mota.
\newblock Factors in randomly perturbed hypergraphs.
\newblock {\em Random Structures \& Algorithms}, 60(2):153--165, 2022.

\bibitem{chang2023powers}
Y.~Chang, J.~Han, and L.~Thoma.
\newblock On powers of tight {H}amilton cycles in randomly perturbed
  hypergraphs.
\newblock {\em Random Structures \& Algorithms}, 63(3):591--609, 2023.

\bibitem{Dirac}
G.~A. Dirac.
\newblock Some theorems on abstract graphs.
\newblock {\em Proc. London Math. Soc.}, 3-2(1):69--81, 1952.

\bibitem{WOS:000287601500007}
A.~Dudek and A.~Frieze.
\newblock Loose {H}amilton cycles in random uniform hypergraphs.
\newblock {\em Electron. J. Combin.}, 18(1), 2011.

\bibitem{dudek2013tight}
A.~Dudek and A.~Frieze.
\newblock Tight {H}amilton cycles in random uniform hypergraphs.
\newblock {\em Random Structures \& Algorithms}, 42(3):374--385, 2013.

\bibitem{MR4052848}
A.~Dudek, C.~Reiher, A.~Ruci\'nski, and M.~Schacht.
\newblock Powers of {H}amiltonian cycles in randomly augmented graphs.
\newblock {\em Random Structures \& Algorithms}, 56(1):122--141, 2020.

\bibitem{MR125031}
P.~{Erd{\H{o}}s and A. R\'enyi}.
\newblock On the evolution of random graphs.
\newblock {\em Magyar Tud. Akad. Mat. Kutat\'o{} Int. K\"ozl.}, 5:17--61, 1960.

\bibitem{erdHos1983supersaturated}
P.~Erd{\H{o}}s and M.~Simonovits.
\newblock Supersaturated graphs and hypergraphs.
\newblock {\em Combinatorica}, 3:181--192, 1983.

\bibitem{MR4298747}
K.~Frankston, J.~Kahn, B.~Narayanan, and J.~Park.
\newblock Thresholds versus fractional expectation-thresholds.
\newblock {\em Ann. of Math. (2)}, 194(2):475--495, 2021.

\bibitem{HPER}
J.~Han, P.~Morris, and A.~Treglown.
\newblock Tilings in randomly perturbed graphs: Bridging the gap between
  {H}ajnal-{S}zemer\'{e}di and {J}ohansson-{K}ahn-{V}u.
\newblock {\em Random Structures \& Algorithms}, 58(3):480--516.

\bibitem{han2020hamiltonicity}
J.~Han and Y.~Zhao.
\newblock Hamiltonicity in randomly perturbed hypergraphs.
\newblock {\em J. Combin. Theory, Ser. B}, 144:14--31, 2020.

\bibitem{MR1782847}
S.~Janson, T.~{\L}uczak, and A.~Rucinski.
\newblock {\em Random graphs}.
\newblock Wiley-Interscience Series in Discrete Mathematics and Optimization.
  Wiley-Interscience, New York, 2000.

\bibitem{factorth}
A.~Johansson, J.~Kahn, and V.~Vu.
\newblock Factors in random graphs.
\newblock {\em Random Structures \& Algorithms}, 33(1):1--28, 2008.

\bibitem{MR4052851}
F.~Joos and J.~Kim.
\newblock Spanning trees in randomly perturbed graphs.
\newblock {\em Random Structures \& Algorithms}, 56(1):169--219, 2020.

\bibitem{MR2312440}
J.~Kahn and G.~Kalai.
\newblock Thresholds and expectation thresholds.
\newblock {\em Combin. Probab. Comput.}, 16(3):495--502, 2007.

\bibitem{MR4273128}
J.~Kahn, B.~Narayanan, and J.~Park.
\newblock The threshold for the square of a {H}amilton cycle.
\newblock {\em Proc. Amer. Math. Soc.}, 149(8):3201--3208, 2021.

\bibitem{KELLY2024507}
T.~Kelly, A.~M\"{u}yesser, and A.~Pokrovskiy.
\newblock Optimal spread for spanning subgraphs of {D}irac hypergraphs.
\newblock {\em J. Combin. Theory, Ser. B}, 169:507--541, 2024.

\bibitem{MR434878}
A.~D. Kor\v{s}unov.
\newblock Solution of a problem of {P}. {E}rd{\H{o}}s and {A}. {R}\'enyi on
  {H}amiltonian cycles in undirected graphs.
\newblock {\em Dokl. Akad. Nauk SSSR}, 228(3):529--532, 1976.

\bibitem{krivelevich2016cycles}
M.~Krivelevich, M.~Kwan, and B.~Sudakov.
\newblock Cycles and matchings in randomly perturbed digraphs and hypergraphs.
\newblock {\em Combin., Probab. and Comput.}, 25(6):909--927, 2016.

\bibitem{MR3595872}
M.~Krivelevich, M.~Kwan, and B.~Sudakov.
\newblock Bounded-degree spanning trees in randomly perturbed graphs.
\newblock {\em SIAM J. Discrete Math.}, 31(1):155--171, 2017.

\bibitem{WOS:000450299200005}
A.~McDowell and R.~Mycroft.
\newblock Hamilton $\ell$-cycles in randomly perturbed hypergraphs.
\newblock {\em Electron. J. Combin.}, 25(4), 2018.

\bibitem{MR3998769}
R.~Montgomery.
\newblock Spanning trees in random graphs.
\newblock {\em Adv. Math.}, 356:106793, 92, 2019.

\bibitem{Nenadov2020}
R.~Nenadov and Y.~Pehova.
\newblock On a {R}amsey-{T}ur\'an variant of the {H}ajnal-{S}zemer\'edi
  theorem.
\newblock {\em SIAM J. Discrete Math.}, 34(2):1001--1010, 2020.

\bibitem{NT}
R.~Nenadov and M.~Truji\'{c}.
\newblock Sprinkling a few random edges doubles the power.
\newblock {\em SIAM J. Discrete Math.}, 35(2):988--1004, 2021.

\bibitem{parczyk2016spanning}
O.~Parczyk and Y.~Person.
\newblock Spanning structures and universality in sparse hypergraphs.
\newblock {\em Random Structures \& Algorithms}, 49(4):819--844, 2016.

\bibitem{MR4654612}
J.~Park and H.~T. Pham.
\newblock A proof of the {K}ahn-{K}alai conjecture.
\newblock {\em J. Amer. Math. Soc.}, 37(1):235--243, 2024.

\bibitem{pham2023toolkit}
H.~T. Pham, A.~Sah, M.~Sawhney, and M.~Simkin.
\newblock A toolkit for robust thresholds.
\newblock {\em arXiv: 2210.03064}, 2023.

\bibitem{MR389666}
L.~P\'osa.
\newblock Hamiltonian circuits in random graphs.
\newblock {\em Discrete Math.}, 14(4):359--364, 1976.

\bibitem{MR1762785}
O.~Riordan.
\newblock Spanning subgraphs of random graphs.
\newblock {\em Combin. Probab. Comput.}, 9(2):125--148, 2000.

\bibitem{MR2195584}
V.~R\"odl, A.~Ruci\'nski, and E.~Szemer\'edi.
\newblock A {D}irac-type theorem for 3-uniform hypergraphs.
\newblock {\em Combin. Probab. Comput.}, 15(1-2):229--251, 2006.

\bibitem{Talagrand2010}
M.~Talagrand.
\newblock Are many small sets explicitly small?
\newblock In {\em S{TOC}'10---{P}roceedings of the 2010 {ACM} {I}nternational
  {S}ymposium on {T}heory of {C}omputing}, pages 13--35. ACM, New York, 2010.

\end{thebibliography}

\appendix
\section{Proof of Lemma~\ref{changabsorbplus}}
In this section, we prove Lemma~\ref{changabsorbplus}.
We first note the following bipartite templete lemma by Montgomery~\cite{MR3998769}, in the following version stated by Nenadov and Pehova~\cite{Nenadov2020}.

\begin{lemma}[\cite{Nenadov2020}, Lemma 2.3]\label{bigraph}
	Let $0<\beta\le 1$ be given.
    There exists $b_0$ such that the following holds for every $b\ge b_0$.
    There exists a bipartite graph $B$ with vertex classes $X_b\cup Y_b$ and $Z_b$ and maximum degree $\Delta(B)\le 40$, such that $|X_b|=b+\beta b$,    $|Y_b|=2b$ and $|Z_b|=3b$, and for every subset $X_b'\subseteq X_b$ with $|X_b'|=b$, the induced graph $B[X_b'\cup Y_b, Z_b]$ contains a perfect matching.
\end{lemma}

We also recall the following classical result on bipartite edge colorings.
\begin{lemma}[\cite{bondy1976graph}, Theorem 6.1]\label{Hall}
	If $G$ is bipartite, then  the  chromatic index of $G$ is equal to its maximum degree $\Delta(G)$.
\end{lemma}

\begin{proof}[{\noindent Proof of Lemma \ref{changabsorbplus}.}]
We choose $\gamma_1$ and $\beta$ such that
$1/n\ll \xi\ll \beta \ll \gamma_1 \ll \tau \ll  \gamma \ll 1/m$.
By property  (i), $\tau n/4\le |V_1|\le \tau n$.
	Let $X=V_1$, $b=|X|/(1+\beta)$ and let $B$ be the bipartite graph given by Lemma  \ref{bigraph} with vertex classes $X_b\cup Y_b$ and $Z_b$. Choose disjoint sets $Y,Z\subseteq V(G^*)\setminus V_1$
	such that $|Y|=2b$ and $|Z|=3b(m-1)$. We now partition $Z$ arbitrarily into $(m-1)$-subsets $\mathcal{Z}=\{Z_i\}_{i\in [3b]}$ and fix bijections $\phi_1: X_b\cup Y_b \rightarrow X\cup Y$ and  $\phi_2: Z_b  \rightarrow \mathcal{Z}$ such that $\phi_1(X_b)=X$ and $\phi_1(Y_b)=Y$.
	
	We first claim that there exists a family $\{A_e\}_{e\in E(B)}$ of pairwise vertex-disjoint $C_1$-subsets of $V(G^*)\setminus (X\cup Y \cup Z)$ such that for every $e=\{w_1, w_2\}\in E(B)$ with $w_1\in X_b\cup Y_b$ and $w_2\in Z_b$, the set $A_e$ is a $(\{\phi_1(w_1)\}\cup\{\phi_2(w_2)\} )$-absorber.  For each edge $e=\{w_1, w_2\} \in E(B)$, let $\phi'(e)$ be the $m$-set $\{\phi_1(w_1)\} \cup \{\phi_2(w_2)\}$ of $V(G^*)$.
    By Lemma \ref{Hall}, $B$ can be decomposed into at most $40$ edge-disjoint matchings.
	We partition each matching into edge-disjoint subsets of size at most  $\gamma_1 n$ and let $M_1, \ldots, M_k$ be the resulting matching, where $k \le   \frac{120\tau }{\gamma_1 }  $.
For each $j \in [k]$, the set $\{\phi'(e) \mid e \in M_j\}$ is a collection of disjoint $m$-sets in $X\cup Y\cup Z$, denoted by $\{S_{j,1},\ldots,S_{j,t_j}\}$, where $t_j \le \gamma_1 n$.
	The families $\{S_{1,1},\ldots,S_{1,t_1}\}, \ldots, \{S_{k,1},\ldots,S_{k,t_k}\}$ are   contained in $X\cup Y\cup Z$, so by the property (iii), w.h.p. there exist vertex-disjoint $C_1$-sets $\{T_{i, j} \mid i \in [k], j \in [t_i]\}$ in  $V(G^*)\setminus (\bigcup_{i\in [k], j\in [t_i]} S_{i,j})$ such that $T_{i, j}$ is an $S_{i, j}$-absorber. By Lemma \ref{bigraph}, we can see there is no isolated vertex in $B$. Thus $X\cup Y\cup Z=\bigcup_{i\in [k], j\in [t_i]} S_{i,j}$ and so $\bigcup_{i\in [k], j\in [t_i]} T_{i,j} \subseteq V(G^*)\setminus (X\cup Y\cup Z)$.
	Then by letting $A_e = T_{i, j}$ where $\phi'(e)=S_{i, j}$, we obtain the desired set of absorbers $\{A_e\}_{e\in E(B)}$.

	Let $V_0=X\cup Y\cup Z\cup (\bigcup_{e\in E(B)} A_e)$.
	We have  $|V_0|\le3m\tau n+120\tau n \cdot C_1 \le \gamma n$ since $\tau \ll \gamma$. We claim that $V_0$ is an absorbing set.  Consider any $R\subseteq V(G^*)\setminus V_0$ with $|R|+|V_0|\in m\mathbb{N}$ and $0\le |R|\le \xi n$,  we first want to  find a $Q\subseteq X$ with $|Q|=\beta b$ and such that w.h.p. $G^*[Q\cup R]$ contains an $F$-factor.

	Note that  the $C_1$-set $T_{i,j}$ is an $S_{i,j}$-absorber, we get $C_1\in m\mathbb{N}$. Thus	   $\beta b+|R| \in m\mathbb{N}$ and  $\beta b-(m-1)|R| \in m \mathbb{N}$.
	Take $Q'\subseteq X$ with $|Q'|=\frac{\beta b-(m-1)|R|}{m}$. Then $|Q'|+|R|\le \frac{\beta b-(m-1)|R|}{m} +\xi n\le \gamma_1 n$ by $\xi \ll \beta \ll \gamma_1$.
	For $Q'\cup R$, by property (ii),  w.h.p. there exists a set $Q''\subseteq X\setminus (Q'\cup R)$ of size $(m-1)|Q'\cup R|$ such that $G^*[Q''\cup Q'\cup R]$ contains an $F$-factor. Let $Q=Q'\cup Q''\subseteq X$. Then $G^*[Q\cup R]$ contains an $F$-factor and $|Q|=|Q'|+(m-1)\cdot (|Q'|+|R|)= \beta b$.

 Let $X'=X\setminus Q$ and $X_b'=\phi_1^{-1}(X')$. Then $|X'|=|X_b'|=b$.
	By Lemma \ref{bigraph}, there is a perfect matching $M$ in $B$ between $X_b'\cup Y_b$ and $Z_b$. For each edge $e\in M$, take an $F$-factor in $G^*[\phi'(e) \cup A_e] $ and for each $e\in E(B)\setminus  M$, take an $F$-factor in $G^*[A_e]$. Therefore, this gives an $F$-factor in $G^*[V_0\setminus Q]$. Thus, w.h.p. $G^*[V_0\cup R]$ contains an $F$-factor.
	This completes the proof.
\end{proof}

\end{document}